\documentclass[final,1p,times]{elsarticle}

\usepackage{amssymb}
 \usepackage{amsthm}
\usepackage{amscd}
\usepackage{amsmath}
\usepackage{amsfonts}
\usepackage{amssymb}
\usepackage{graphicx}

\numberwithin{equation}{section}

\newcommand{\ra}{\rightarrow}
\newcommand{\p}{\partial}
\newcommand{\f}{\frac}

\newcommand{\be}{\begin{equation}}
\renewcommand{\ra}{\rightarrow}
\newcommand{\ee}{\end{equation}}
\newcommand{\bea}{\begin{eqnarray}}
\newcommand{\eea}{\end{eqnarray}}
\newcommand{\bna}{\begin{eqnarray*}}
\newcommand{\ena}{\end{eqnarray*}}

\renewcommand{\le}{\left}
\newcommand{\ri}{\right}

\journal{***}

\begin{document}

\begin{frontmatter}

\title{Trudinger-Moser inequalities on complete noncompact Riemannian manifolds}

\author{Yunyan Yang}
 \ead{yunyanyang@ruc.edu.cn}
\address{ Department of Mathematics,
Renmin University of China, Beijing 100872, P. R. China}

\begin{abstract}
 Let $(M,g)$ be a complete noncompact Riemannian $n$-manifold ($n\geq 2$).
 If there exist positive constants $\alpha$, $\tau$ and $\beta$ such
 that
 $$\sup_{u\in W^{1,n}(M),\,\|u\|_{1,\tau}\leq 1}
 \int_M\le(e^{\alpha|u|^{\f{n}{n-1}}}-\sum_{k=0}^{n-2}\f{\alpha^k|u|^{\f{nk}{n-1}}}{k!}\ri)dv_g
 \leq \beta,$$
 where $\|u\|_{1,\tau}=\|\nabla_gu\|_{L^n(M)}+\tau\|u\|_{L^n(M)}$,
 then we say that Trudinger-Moser inequality holds.
 Suppose Trudinger-Moser inequality holds, we prove that
 there exists some positive constant
 $\epsilon$ such that ${\rm Vol}_g(B_x(1))\geq \epsilon$ for all $x\in M$. Also we give a sufficient
 condition under which Trudinger-Moser inequality holds, say the Ricci curvature of $(M,g)$
 has lower bound and its injectivity radius is positive. Moreover, Adams inequality is
 discussed in this paper. For application of Trudinger-Moser inequalities, we obtain existence
 results for some quasilinear equations with nonlinearity of exponential growth.

\end{abstract}


\end{frontmatter}

\section{Introduction}

Let $\Omega$ be a smooth bounded domain in $\mathbb{R}^n$ ($n\geq
2)$ and $C_0^\infty(\Omega)$ be a space of smooth functions with
compact support in $\Omega$. Let $W_0^{m,p}(\Omega)$ be the
completion of $C_0^\infty(\Omega)$ under the Sobolev norm
\be\label{S-n}
 \|u\|_{W_0^{m,p}(\Omega)}:=\le(\sum_{l=0}^m\int_\Omega|\nabla^lu|^pdx\ri)^{1/p}.
\ee Assume that $m$ is an integer satisfying $1\leq m<n$. Then
Sobolev embedding theorem asserts that
$W_0^{m,p}(\Omega)\hookrightarrow L^q(\Omega)$, $1\leq q\leq
{np}/{(n-mp)}$. Concerning the limiting case $mp=n$, one has
$W_0^{m,n/m}(\Omega)\hookrightarrow L^q(\Omega)$ for all $q\geq 1$.
But the embedding is not valid for $q=\infty$. To fill this gap, it
is natural to find the maximal growth function
$g:\mathbb{R}\ra\mathbb{R}^+$ such that
$$\sup_{u\in W_0^{m,n/m}(\Omega),\,\|u\|_{W_0^{m,n/m}(\Omega)}\leq 1}\int_\Omega g(u)dx<\infty.$$
In the case $m=1$, Trudinger \cite{Trudinger} and Pohozaev \cite{Pohozaev}
found independently that the maximal growth is of exponential type.
More precisely, there exist two positive constants $\alpha_0$ and
$C$ depending only on $n$ such that \be\label{T-Mos}\sup_{u\in
W_0^{1,n}(\Omega),\,\|u\|_{W_0^{1,n}(\Omega)}\leq 1}\int_\Omega
e^{\alpha_0|u|^{\f{n}{n-1}}}dx\leq C|\Omega|,\ee where $|\Omega|$
denotes the Lebesgue measure of $\Omega$.  Moser \cite{Moser}
obtained the best constant $\alpha_n=n\omega_{n-1}^{1/(n-1)}$ such
that the above supremum is finite when $\alpha_0$ is replaced by
$\alpha_n$, where $\omega_{n-1}$ is the area of the unit sphere in
$\mathbb{R}^n$.  Moser's work relies on a rearrangement argument
\cite{Hardy}. In literature the kind of inequalities like
(\ref{T-Mos}) are called Trudinger-Moser inequalities.

  Adams \cite{Adams}
generalized inequality (\ref{T-Mos}) to the case of general $m:
1\leq m<n$ as follows. For any $u\in W_0^{m,n/m}(\Omega)$, the
$l$-th order gradient of $u$ reads \be\label{gradient}
 \nabla^lu=\le\{\begin{array}{lll}
  \Delta^{\f{l}{2}}u,\quad&{\rm if}\,\,\,l\,\,\,{\rm is}\,\,\, {\rm
  even,}\\[1.5ex]\nabla\Delta^{\f{l-1}{2}}u,&{\rm if}\,\,\,l\,\,\,{\rm is}\,\,\, {\rm
  odd,}\end{array}\ri.
\ee there exits a positive constant $C_{m,n}$ such that
\be\label{Ada}
 \sup_{u\in W_0^{m,n/m}(\Omega),\,\|u\|_{W_0^{m,n/m}(\Omega)}\leq 1}\int_\Omega
 e^{\beta_0|u|^{\f{n}{n-m}}}dx\leq C_{m,n}|\Omega|,
\ee where $\beta_0$  is the best constant depending only on $n$ and
$m$, namely
\be\label{beta}\beta_0=\beta_0(m.n):=\le\{\begin{array}{lll}\f{n}{\omega_{n-1}}\le[\f{\pi^{{n}/{2}}2^m\Gamma\le({(m+1)}/{2}\ri)}
{\Gamma\le({(n-m+1)}/{2}\ri)}\ri]\quad&{\rm when}\,\,\,m\,\,\,{\rm
is}\,\,\,odd
\\[1.5ex]\f{n}{\omega_{n-1}}\le[\f{\pi^{{n}/{2}}2^m\Gamma\le({m}/{2}\ri)}
{\Gamma\le({(n-m)}/{2}\ri)}\ri]&{\rm when}\,\,\,m\,\,\,{\rm
is}\,\,\,even.\end{array}\ri.\ee
  The inequality (\ref{Ada}) is known as Adams inequality.  Adams first represented a function $u$ in terms of its gradient function
$\nabla^mu$ by using a convolution operator. Then using the O'Neil's
idea \cite{Neil} of rearrangement of convolution of two functions
and the idea which originally goes back to Garcia, he obtained (\ref{Ada}).

There are many types of extensions for Trudinger-Moser inequality
and Adams inequality.  One is to establish such inequalities
 on the whole euclidian space $\mathbb{R}^n$.
 Cao \cite{Cao} employed the decreasing rearrangement argument to
 prove that for all $\alpha<4\pi$ and $A>0$, there
exists a constant $C$ depending only on $\alpha$ and $A$ such that
for all $u\in W^{1,2}(\mathbb{R}^2)$ with
$\int_{\mathbb{R}^2}|\nabla u|^2dx\leq
1,\,\int_{\mathbb{R}^2}u^2dx\leq A$, there holds
\be\label{cao}\int_{\mathbb{R}^2}\le(e^{\alpha u^2}-1\ri)dx\leq
C.\ee
 His argument was generalized
 to $n$-dimensional case by
 do \'O \cite{doo1} and  Panda \cite{Panda} independently. Later,
 Adachi-Tanaka \cite{AT} gave another
type of generalization. All these inequalities are subcritical ones
since $\alpha<\alpha_n$. It was Ruf \cite{Ruf} who first proved the
critical Trudinger-Moser inequality in the whole euclidian space
$\mathbb{R}^2$ and gave out extremal functions via more delicate
analysis.  This result was generalized to $n$-dimensional case by
Li-Ruf \cite{LiRuf} through combining symmetrization and blow-up
analysis. Subsequently, using the decreasing rearrangement argument
and Young's inequality,
Adimurthi-Yang \cite{Adi-Yang} derived an
interpolation of Trudinger-Moser inequality and Hardy inequality in
$\mathbb{R}^n$, which can be viewed as a singular Trudinger-Moser
inequality. Another kind of singular Trudinger-Moser inequality was
recently established by Wang-Ye \cite{Wang-Ye} through the method of
blow-up analysis.

Substantial progresses on Adams inequality in $\mathbb{R}^n$ was
also made recently. Following lines of Adams, Kozono et al.
\cite{Kozono} obtained subcritical Adams inequality in the whole
euclidian space $\mathbb{R}^n$. Based on rearrangement argument of
Trombetti-Vazquez \cite{T-V}, Ruf-Sani \cite{Ruf-Sani} proved the
critical Adams inequalities in $\mathbb{R}^n$ as follows. Let $m$ be
an even integer less than $n$. Assume that $u\in
W_0^{m,n/m}(\mathbb{R}^n)$ and
$\|(-\Delta+I)^{m/2}u\|_{L^{n/m}(\mathbb{R}^n)}\leq 1$. There exists
a constant $C>0$ depending only on $n$ and $m$ such that
$$\int_{\mathbb{R}^n}
\le(e^{\beta_0|u|^{\f{n}{n-m}}}-\sum_{k=0}^{j-2}\f{\beta_0^k|u|^{\f{nk}{n-m}}}{k!}\ri)dx<C,$$
where $j$ is the smallest integer great than or equal to $n/m$.

Another extension is to establish Trudinger-Moser inequality and
Adams inequality on compact Riemannian manifolds. Let $(M,g)$ be a
compact Riemannian $n$-manifold. For $u\in W^{1,n}(M)$, it was shown
by Aubin \cite{Aubin10} that
$\exp(\alpha|u|^{n/(n-1)}\|u\|_{W^{1,n}(M)}^{-n/(n-1)})$ is
integrable for sufficiently small $\alpha>0$ which does not depend
on $u$. In fact, this is an easy consequence of Trudinger-Moser
inequality and finite partition of unity on $M$. Let
$\tilde{\alpha}$ be the supremum of the above $\alpha$'s. It was
first found by Cherrier \cite{cherrier1} that
$\tilde{\alpha}=\alpha_n$. Cherrier \cite{cherrier2} obtained
similar results for $u\in W_{m,n/m}(M)$. Following the lines of
Adams, Fontana \cite{Fontana} obtained critical Adams inequality on $(M,g)$.
In 1997, using the method of blow-up analysis, Ding et al.
\cite{DJLW} established a nice Trudinger-Moser inequality on compact
Riemannian surface and successfully applied it to deal with the
prescribed Gaussian curvature problem. Adapting the argument of Ding
et al., Li \cite{LiJPDE,LiSci}  and Li-Liu \cite{Li-Liu} proved the
existence of extremal functions for Trudinger-Moser inequalities.
Their idea was also employed by the author
\cite{YangIJM,YangTrans,YangMZ} to find extremal functions for
various Trudinger-Moser type inequalities. For vector bundles over a
compact Riemannian 2-manifold, Li-Liu-Yang obtained Trudinger-Moser
inequalities in \cite{Li-Liu-Yang}.

Among other contributions, we mention the following results. Using
the method of blow-up analysis, Adimurthi-Druet \cite{A-D} proved
that when $0\leq\alpha<\lambda_1(\Omega)$, there holds
$$\sup_{u\in W_0^{1,2}(\Omega),\,\|\nabla u\|_{2}\leq 1}\int_\Omega
e^{4\pi u^2(1+\alpha\|u\|_{2}^2)}dx<\infty,$$
 where $\lambda_1(\Omega)$ is
the first eigenvalue of Laplacian on bounded smooth domain
$\Omega\subset\mathbb{R}^2$.
 Moreover, the supremum
is infinite when $\alpha\geq \lambda_1(\Omega)$. Later this result
was generalized by the author \cite{YangJFA} and Lu-Yang
\cite{Lu-Yang,Lu-Yang1,Lu-Yang2}. 

 Although there are fruitful results on euclidian space and compact Riemannian manifolds,
 we know little about Trudinger-Moser
 inequalities on complete noncompact Riemannian manifolds.
 In this paper, we concern this problem. Let $(M,g)$ be any complete noncompact
Riemannian $n$-manifold. Throughout this paper, all the manifolds
are assumed to be without boundary, and of
dimension $n\geq 2$. We say that Trudinger-Moser inequality holds on
$(M,g)$ if there exist positive constants $\alpha$, $\tau$ and
$\beta$ such that
 \be\label{Tm-hold}\sup_{u\in W^{1,n}(M),\,\|u\|_{1,\tau}\leq 1}
 \int_M\le(e^{\alpha|u|^{\f{n}{n-1}}}-\sum_{k=0}^{n-2}\f{\alpha^k|u|^{\f{nk}{n-1}}}{k!}\ri)dv_g
 \leq \beta,\ee
 where \be\label{1tau}\|u\|_{1,\tau}=\le(\int_M|\nabla_g u|^ndv_g\ri)^{1/n}+
 \tau\le(\int_M|u|^ndv_g\ri)^{1/n}.\ee
 If the above
 supremum is infinite for all $\alpha>0$ and $\tau>0$, then we say that Trudinger-Moser
 inequality is not valid on $(M,g)$.
 Motivated by Sobolev embedding (Hebey \cite{Hebey}, Chapter
3), in this paper, we propose and answer the
  following three questions.\\

  \noindent{ $(Q_1)$} {\it Which kind of complete noncompact Riemannian manifolds can possibly make
  Trudinger-Moser inequalities hold?}

  \noindent{ $(Q_2)$} {\it What geometric assumptions should we consider in order to obtain
  Trudinger-Moser inequalities on complete noncompact Riemannian manifolds?}

  \noindent{ $(Q_3)$} {\it Are those geometric assumptions in $(Q_2)$ necessary?}\\

   This paper is organized as follows: In Section 2, we state our main results.
   From section 3 to section 5, we answer the questions
   $(Q_1)$-$(Q_3)$,
  respectively. Adams inequalities are considered in section 6. Finally, Trudinger-Moser inequalities
  are applied to nonlinear analysis in section 7.

  \section{Main results}

   In this section, we answer questions
  $(Q_1)$-$(Q_3)$, and give an application of Trudinger-Moser inequality. Throughout this paper, we denote for
  simplicity a function $\zeta:\mathbb{N}\times
  [0,\infty)\ra\mathbb{R}$ by
  \be\label{zeta}\zeta(l,t)=e^t-\sum_{k=0}^{l-2}\f{t^k}{k!},\quad \forall l\geq 2.\ee
  From (\cite{Yang}, lemma 2.1 and lemma 2.2), we know that
  \be\label{Y1}\le(\zeta(l,t)\ri)^q\leq \zeta(l,qt)\ee
  and
  \be\label{Y2}\zeta(l,t)\leq \f{1}{\mu}\zeta(l,\mu t)+\f{1}{\nu}\zeta(l,\nu t).\ee
  for all $l\geq 2$, $q\geq 1$, $t\in[0,\infty)$, and $\mu>0$, $\nu>0$
  satisfying $1/\mu+1/\nu=1$.\\

  The following proposition answers question $(Q_1)$.\\

\noindent{\bf Proposition 2.1.} {\it Let $(M,g)$ be a complete
Riemannian $n$-manifold. Suppose that Trudinger-Moser inequality
holds on $(M,g)$, i.e. there exist positive constants $\alpha$, $\tau$ and
$\beta$ such that (\ref{Tm-hold}) holds. Then the Sobolev space
$W^{1,n}(M)$ is embedded in $L^q(M)$ continuously for any $q\geq n$.
 Furthermore, for any $r>0$ there exists a positive constant $\epsilon$
 depending only on $n$, $\alpha$, $\tau$, $\beta$ and $r$
 such that ${\rm Vol}_g(B_x(r))\geq \epsilon$ for all $x\in M$, where $B_x(r)$ denotes the geodesic ball
 centered at $x$ with radius $r$.}\\

 From proposition 2.1 we know that there are indeed complete noncompact Riemannian manifolds such that
  Trudinger-Moser inequalities are not valid, namely \\

 \noindent{\bf Corollary 2.2.} {\it For any integer $n\geq 2$, there exists a complete noncompact Riemannian
 $n$-manifold on which Trudinger-Moser inequality is not valid.}\\

  To answer question $(Q_2)$, we have the following:\\

 \noindent{\bf Theorem 2.3.} {\it Let $(M,g)$ be a complete noncompact Riemannian $n$-manifold.
 Suppose that its Ricci curvature has lower bound, namely
 ${\rm Rc}_{(M,g)}\geq Kg$ for some constant $K\in\mathbb{R}$, and its injectivity radius is strictly positive,
 namely ${\rm inj}_{(M,g)}\geq i_0$ for some constant $i_0>0$. Then
 we have

 \noindent $(i)$ for any $0\leq\alpha<\alpha_n=n\omega_{n-1}^{1/(n-1)}$, there
  exists positive constants $\tau$ and $\beta$ depending only on $n$, $\alpha$, $K$ and $i_0$ such that
  (\ref{Tm-hold}) holds.
 As a consequence, $W^{1,n}(M)$ is embedded in $L^q(M)$ continuously
 for any $q\geq n$;

 \noindent $(ii)$ for any $\alpha>\alpha_n$ and any $\tau>0$, the supremum in (\ref{Tm-hold}) is
 infinite;

 \noindent $(iii)$ for any $\alpha>0$ and any $u\in W^{1,n}(M)$,
 there holds $\zeta(n,\alpha|u|^{n/(n-1)})\in L^1(M)$.}\\

Now we turn to question $(Q_3)$. The following proposition
implies that one of the hypotheses of theorem 2.3, the injectivity
radius is strictly positive, can
not be removed. \\

\noindent{\bf Proposition 2.4.} {\it For any integer $n\geq 2$,
there exists a complete noncompact Riemannian
 $n$-manifold, whose Ricci curvature has lower bound, such that Trudinger-Moser inequality is not
 valid on it.}\\

  We shall construct complete noncompact Riemannian manifolds on which
 Trudinger-Moser inequalities hold, but their Ricci curvatures are unbounded from below.
 This implies that the other hypothesis of theorem 2.3, Ricci curvature has lower bound,
 is not necessary. Namely\\

 \noindent{\bf Proposition 2.5.} {\it For any integer $n\geq 2$, there
exists a complete noncompact Riemannian
 $n$-manifold on which
 Trudinger-Moser inequality holds, but its Ricci curvature is unbounded from
 below.}\\

 Let us explain the idea of proving proposition 2.1 and theorem 2.3.
 The first part of conclusions of proposition 2.1, $W^{1,n}(M)\hookrightarrow L^q(M)$ for all $q\geq n$,
  is based on an observation
  $$\int_M\zeta\le(n,\alpha|u|^{\f{n}{n-1}}\ri)dv_g=\sum_{k=n-1}^\infty\f{\alpha^k}{k!}\int_M|u|^{\f{nk}{n-1}}dv_g.$$
 To find some $\epsilon>0$ such that ${\rm Vol}_g(B_x(r))\geq \epsilon$ for all
 $x\in M$, we employ the method of Carron (\cite{Hebey}, lemma 3.2) who obtained similar
 result for Sobolev embedding. For the proof of theorem 2.3, we first derive a uniform local Trudinger-Moser
 inequality (lemma 4.2 below). Then using harmonic coordinates and
 Gromov's covering lemma, we get the desired global Trudinger-Moser
 inequality. The proofs of corollary 2.2, proposition 2.4 and proposition 2.5 are all based on construction
 of Riemannian manifolds. \\

 Concerning Adams inequalities on complete noncompact Riemannian
 manifolds, we have the following:\\

 \noindent{\bf Theorem 2.6.} {\it Let $(M,g)$ be a complete noncompact Riemannian $n$-manifold.
 Suppose that there exist positive constants $C(k)$ and $i_0$  such that
 $|\nabla_g^{k}{\rm Rc}_{(M,g)}|\leq C(k)$, $k=0,1,\cdots,m-1$,  ${\rm inj}_{(M,g)}\geq
 i_0>0$. Let $j=n/m$ when
 $n/m$ is an integer, and $j=[n/m]+1$ when $n/m$ is not an
 integer, where $[n/m]$ denotes the integer part of $n/m$.
 Then we conclude the following:\\
 \noindent$(i)$ there exist positive constants $\alpha_0$ and $\beta$ depending only on $n$, $m$,
 $C(k)$, $k=1,\cdots,m-1$, and $i_0$
 such that
 $$\sup_{\|u\|_{W^{m,n/m}(M)}\leq 1}\int_M\zeta\le(j,\alpha_0|u|^{\f{n}{n-m}}\ri)dv_g\leq \beta.$$
 As a consequence, $W^{m,{n}/{m}}(M)$ is embedded in $L^q(M)$ continuously
 for any $q\geq {n}/{m}$;\\
 \noindent$(ii)$ for any $\alpha>0$ and any $u\in W^{m,n/m}(M)$,
 there holds $\zeta(j,\alpha|u|^{n/(n-m)})\in L^1(M)$.}\\

 The proof of theorem 2.6 is similar to that of theorem 2.3. It should be remarked that the
 existing proofs of Trudinger-Moser
 inequalities or Adams inequalities for the euclidian space $\mathbb{R}^n$
 are all based on rearrangement argument, which is difficult to be applied to complete
 noncompact Riemannian manifold case. Our method is from uniform local estimates to global
 estimates. It does not depend on the rearrangement theory directly.\\

 Trudinger-Moser inequality plays an important role in nonlinear analysis. Let $(M,g)$
 be a complete noncompact Riemannian $n$-manifold. $\nabla_g$ denotes its covariant derivative, and
 ${\rm div}_g$ denotes its divergence operator.  Assume the Ricci
 curvature of $(M,g)$  has lower bound and the injectivity radius is
 strictly positive. We consider the existence results for
 the following quasilinear equation.
 \be\label{equa}-{\rm div}_g(|\nabla_g u|^{n-2}\nabla_g u)+v(x)|u|^{n-2}u=\phi(x)f(x,u),\ee
 where $v(x)$, $\phi(x)$ and $f(x,t)$ are all continuous functions, and $f(x,t)$ behaves like $e^{\alpha t^{n/(n-1)}}$ as
 $t\ra+\infty$. In the case that $(M,g)$ is the standard euclidean space
 $\mathbb{R}^n$ and $\phi(x)=|x|^{-\beta}$ $(0\leq\beta<n)$, problem (\ref{equa}) has been studied by do
 \'O et. al. \cite{doo,do-de}, Adimurthi-Yang \cite{Adi-Yang},
 Yang \cite{Yang}, Lam-Lu \cite{LamLu} and Zhao \cite{Zhao}. Let $O$ be a fixed point
 of $M$ and $d_g(\cdot,\cdot)$ be the geodesic distance between two
 points of $(M,g)$.
 Assume that $\phi(x)$ satisfies the
 following hypotheses.\\

 \noindent{$(\phi_1)$ $\phi(x)\in L^p_{\rm loc}(M)$ for some $p>1$, i. e., for any $R>0$
 there holds $\phi(x)\in L^p(B_O(R))$;}\\
 $(\phi_2)$ $\phi(x)> 0$ for all $x\in M$ and there exist positive constants $C_0$ and $R_0$
 such that $\phi(x)\leq C_0$ for all $x\in M\setminus B_O(R_0)$.\\

  The potential $v(x)$ is assumed to satisfy the following:\\

 \noindent $(v_1)$ there exists some constant $v_0>0$ such that $v(x)\geq
 v_0$ for all $x\in M$;\\
 $(v_2)$ either $v(x)\in L^{1/{(n-1)}}(M)$ or $v(x)\ra+\infty$ as
 $d_g(O,x)\ra+\infty$.\\

 The nonlinearity $f(x,t)$ satisfies the following hypotheses.\\

 \noindent $(f_1)$ there exist constants $\alpha_0$, $b_1$, $b_2>0$
such that for all $(x,t)\in M\times\mathbb{R}^+$,
$$|f(x,t)|\leq
b_1 t^{n-1}+b_2\zeta\le(n,\alpha_0t^{n/(n-1)}\ri);$$
$(f_2)$ there
exists some constant $\mu>n$ such that for all $x\in M$ and $t>0$,
$$0<\mu F(x,t)\equiv \mu\int_0^tf(x,s)ds\leq tf(x,t);$$
 $(f_3)$ there exist constants $R_1$, $A_1>0$ such that if $t\geq
R_1$, then for all $x\in M$ there holds
$$F(x,t)\leq A_1f(x,t).$$

Define a function space
\be\label{space}E=\le\{u\in
W^{1,n}(M):\int_{M}v(x)|u|^ndv_g<\infty\ri\}.\ee We say that $u\in
E$ is a weak solution of problem (\ref{equa}) if for all $\varphi\in
E$ we have
$$\int_{M}\le(|\nabla_g u|^{n-2}\nabla_g u\nabla_g
\varphi+v(x)|u|^{n-2}u\varphi\ri)dv_g
=\int_{M}\phi(x){f(x,u)}\varphi dv_g.$$ Define a weighted eigenvalue
for the $n$-Laplace operator by
\be\label{lamda}\lambda_\phi=\inf_{u\in
E,\,u\not\equiv0}\f{\int_M(|\nabla_g u|^n+v(x)|u|^n)dv_g}
{\int_{M}\phi(x)|u|^ndv_g}.\ee
Then we state the following:\\

 \noindent{\bf Theorem 2.7.} {\it Let $(M,g)$ be a complete noncompact Riemannian $n$-manifold.
 Suppose that ${\rm Rc}_{(M,g)}\geq Kg$ for some constant $K\in\mathbb{R}$, and ${\rm inj}_{(M,g)}\geq i_0$
 for some positive constant $i_0$. Assume that $v(x)$ is a continuous function satisfying $(v_1)$ and
 $(v_2)$,  $\phi(x)$ is a continuous function satisfying $(\phi_1)$ and
 $(\phi_2)$,
 $f:M\times\mathbb{R}\ra\mathbb{R}$ is a continuous
 function and the hypotheses $(f_1)$, $(f_2)$ and $(f_3)$ are satisfied. Furthermore we
 assume

 \noindent $(f_4)$ $\limsup_{t\ra
0+}{nF(x,t)}/{t^n}<\lambda_\phi$ uniformly in $x\in M$;\\
 $(f_5)$ there exist constants $q>n$ and $C_q$ such that for all $(x,t)\in M\times[0,\infty)$
 $$f(x,t)\geq C_q t^{q-1},$$
 where
 $$C_q>\le(\f{q-n}{q}\ri)^{{(q-n)}/{n}}\le(\f{p\alpha_0}{(p-1)\alpha_n}\ri)^{(q-n)(n-1)/n}S_q^q$$
 and
 \be\label{Sq}S_q=\inf_{u\in E\setminus\{0\}}\f{\le(\int_M(|\nabla_g u|^n+v(x)|u|^n)dv_g\ri)^{1/n}}
 {\le(\int_M \phi(x) |u|^qdv_g\ri)^{1/q}}.\ee
 Then the problem (\ref{equa}) has a nontrivial nonnegative weak
 solution.}\\

 \noindent{\bf Remark 2.8.} We shall prove that $S_q$ can be attained (lemma 7.2
 below). When $(M,g)$ is the standard euclidian space $\mathbb{R}^n$, $\phi(x)=|x|^{-\beta}$ for
 $0\leq\beta<n$, $(f_1)$-$(f_4)$ and
 $$\leqno(H_5) \,\, \liminf_{s\ra
+\infty}sf(x,s)e^{-\alpha_0s^\f{n}{n-1}}=\beta_0>\mathcal{M}$$
uniformly in $x$, where $\mathcal{M}$ is some sufficiently large
number, we obtained similar existence result in \cite{Yang}.
The following proposition implies that the set of functions
satisfying $(f_1)$-$(f_5)$ is not empty and assumptions
$(f_1)$-$(f_5)$ do not imply $(H_5)$.\\

 \noindent{\bf Proposition 2.9.} {\it There exist continuous
 functions
 $f:M\times\mathbb{R}\ra\mathbb{R}$ such that
 $(f_1)$-$(f_5)$ are satisfied, but $(H_5)$ is not satisfied.}\\

We also consider multiplicity results for a perturbation of the
problem (\ref{equa}), namely \be\label{equa1}-{\rm div}_g(|\nabla_g
u|^{n-2}\nabla_g u)+v(x)|u|^{n-2}u=\phi(x)f(x,u)+\epsilon h(x),\ee
 where $h(x)\in E^*$, the dual space of $E$. If $h\not\equiv 0$ and $\epsilon>0$ is
 sufficiently small, under some assumptions there exist at least two distinct weak
 solutions to (\ref{equa1}). Precisely, we have the following theorem.\\

 \noindent {\bf Theorem 2.10.} {\it Let $(M,g)$ be a complete noncompact Riemannian $n$-manifold.
 Suppose that ${\rm Rc}_{(M,g)}\geq Kg$ for some constant $K\in\mathbb{R}$, and ${\rm inj}_{(M,g)}\geq i_0$
 for some positive constant $i_0$.  Assume $f(x,t)$ is
continuous in $M\times\mathbb{R}$ and $(f_1)$-$(f_5)$ are satisfied.
Both $v(x)$ and $\phi(x)$ are continuous in $M$ and $(v_1)$,
$(v_2)$, $(\phi_1)$, $(\phi_2)$ are satisfied, $h$ belongs to $E^*$,
the dual space of $E$, with $h\geq 0$ and $h\not\equiv 0$. Then
there exists $\epsilon_0>0$ such that if $0<\epsilon<\epsilon_0$,
then the problem (\ref{equa1}) has at least two distinct nonnegative
weak
solutions.}\\

The proofs of theorem 2.7 and theorem 2.10 are based on theorem 2.3,
Mountain-pass theorem and Ekeland's variational principle. Though
similar idea was used in the case $(M,g)$ is the standard euclidian
space $\mathbb{R}^n$ \cite{Adi-Yang,doo,do-de,LamLu,Yang},
technical difficulties caused by manifold structure must be
smoothed.

\section{Necessary conditions}

In this section, we consider the necessary conditions under
which Trudinger-Moser inequality holds. Precisely we shall prove proposition 2.1 and corollary 2.2.
Firstly we have the following:\\

\noindent{\bf Lemma 3.1.} {\it Let $(M,g)$ be a complete Riemannian
$n$-manifold. Suppose that there exist constants $q>n$, $A>0$ and
$\tau>0$ such that for all $u\in W^{1,n}(M)$, there holds
 \be\label{em1}
 \le(\int_M |u|^qdv_g\ri)^{1/q}\leq A\|u\|_{1,\tau},
\ee where $\|u\|_{1,\tau}$ is defined by (\ref{1tau}).
 Then for any $r>0$ there exists some positive constant $\epsilon$ depending only on
 $A$, $n$, $q$, $\tau$, and $r$ such that
 for all $x\in M$, ${\rm Vol}_g(B_x(r))\geq \epsilon$.}\\

 \noindent{\it Proof.}  Let $r>0$, $x\in M$, and $\phi\in
 W^{1,n}(M)$ be such that $\phi=0$ on $M\setminus B_x(r)$.
 By H\"older's inequality,
 $$\le(\int_M |\phi|^ndv_g\ri)^{1/n}\leq {\rm Vol}_g(B_x(r))^{\f{1}{n}-\f{1}{q}}\le(\int_M |\phi|^qdv_g\ri)^{1/q}.$$
 This together with (\ref{em1}) gives
 \be\label{em3}
 \le(1-\tau A{\rm Vol}_g(B_x(r))^{\f{1}{n}-\f{1}{q}}\ri)\le(\int_M
 |\phi|^qdv_g\ri)^{1/q}\leq A\le(\int_M(|\nabla
 \phi|^ndv_g\ri)^{1/n}.
 \ee
 Fix $x\in M$ and $R>0$. Then either
 \be\label{p1}{\rm Vol}_g(B_x(R))>\le(\f{1}{2\tau A}\ri)^{{nq}/{(q-n)}}\ee
 or
 \be\label{p2}{\rm Vol}_g(B_x(R))\leq\le(\f{1}{2\tau A}\ri)^{{nq}/{(q-n)}}.\ee
 If (\ref{p2}) holds, then we have
 $$1-\tau A{\rm Vol}_g(B_x(R))^{\f{1}{n}-\f{1}{q}}\geq {1}/{2},$$
 and whence for all $r\in(0,R]$ and all $\phi\in W^{1,n}(M)$ with
 $\phi=0$ on $M\setminus B_x(r)$,
 \be\label{p4}\le(\int_M
 |\phi|^qdv_g\ri)^{1/q}\leq 2A\le(\int_M(|\nabla
 \phi|^ndv_g\ri)^{1/n}.\ee
 Now we set
 $$\phi(y)=\le\{\begin{array}{lll}
 r-d_g(x,y)\quad&{\rm when} \quad d_g(x,y)\leq r\\[1.5ex]
 0&{\rm when} \quad d_g(x,y)> r.
 \end{array}\ri.$$
 Clearly $\phi\in W^{1,n}(M)$, $\phi=0$ on $M\setminus B_x(r)$,
 $\phi\geq r/2$ on $B_x(r/2)$,
 and $|\nabla\phi|=1$ almost everywhere in $B_x(r)$.
 It then follows from (\ref{p4}) that
 $$\f{r}{2}{\rm Vol}_g(B_x(r/2))^{1/q}\leq 2A{\rm Vol}_g(B_x(r))^{1/n}.$$
 Hence we have for all $r\leq R$,
 $${\rm Vol}_g(B_x(r))\geq \le(\f{r}{4A}\ri)^n{\rm Vol}_g(B_x(r/2))^{n/q}.$$
 By induction we obtain for any positive integer $m$,
 \be\label{ind}
  {\rm Vol}_g(B_x(R))\geq
  \le(\f{R}{2A}\ri)^{n\alpha(m)}\le(\f{1}{2}\ri)^{n\beta(m)}{\rm
  Vol}_g(B_x(R/2^m))^{(n/q)^m},
 \ee
 where
 $$\alpha(m)=\sum_{j=1}^m({n}/{q})^{j-1},\quad \beta(m)=\sum_{j=1}^m
 j(n/q)^{j-1}.$$
 On one hand we know from (\cite{GaHL}, Theorem 3.98) that
 ${\rm Vol}_g(B_x(r))=\f{\omega_{n-1}}{n}r^n(1+o(r))$,
 where $\omega_{n-1}$ is the area of the euclidean unit sphere in $\mathbb{R}^n$, and
 $o(r)\ra 0$ as $r\ra 0$. One can see without any difficulty that
 $$\lim_{m\ra\infty}{\rm
  Vol}_g(B_x(R/2^m))^{(n/q)^m}=1.$$
 On the other hand we have
 $$\sum_{j=1}^\infty({n}/{q})^{j-1}=\f{q}{q-n},\quad
 \sum_{j=1}^\infty j(n/q)^{j-1}=\f{q^2}{(q-n)^2}.$$
 Hence, passing to the limit $m\ra\infty$ in (\ref{ind}), one
 concludes that
 $${\rm
 Vol}_g(B_x(R))\geq\le(\f{R}{2^{(2q-n)/(q-n)}A}\ri)^{nq/(q-n)}.$$
 This together with (\ref{p1}), (\ref{p2}) implies that
 $${\rm
 Vol}_g(B_x(R))\geq\min\le\{\f{1}{2\tau A},\f{R}{2^{(2q-n)/(q-n)}A}\ri\}^{nq/(q-n)}$$
 and completes the proof of the lemma. $\hfill\Box$\\

 It should be pointed out that the above argument is a modification of that of Carron (\cite{Hebey}, lemma 3.2).
 Note that the condition (\ref{em1}) implies that
 $W^{1,n}(M)$ is continuously embedded in $L^q(M)$ for some $q>n$. This is different from the assumption of (\cite{Hebey}, lemma 3.2).\\

 To prove proposition 2.1, we also need the following interpolation inequality.\\

 \noindent{\bf Lemma 3.2.} {\it Let $\tau$ be any positive real number. Suppose there exist positive constants
 $q_1$, $q_2$, $A_1$ and $A_2$ such that $q_2>q_1>0$ and
 \be\label{cn}\le(\int_M |u|^{q_i}dv_g\ri)^{1/{q_i}}\leq A_i\|u\|_{1,\tau}\ee
 for all $u\in W^{1,n}(M)$, $i=1,2$.
 Then for all $q: q_1<q<q_2$ there exists a positive constant $A=A(A_1,A_2,q_1,q_2)$
 such that
 \be\label{concl}\le(\int_M |u|^{q}dv_g\ri)^{1/{q}}\leq A\|u\|_{1,\tau}\ee
 for all $u\in W^{1,n}(M)$.}\\

 \noindent{\it Proof.} For any $u\in W^{1,n}(M)\setminus\{0\}$, we set
 $\widetilde{u}=u/\|u\|_{1,\tau}$. It follows from (\ref{cn})
 that
 $$\le(\int_M |\widetilde{u}|^{q_i}dv_g\ri)^{1/{q_i}}\leq A_i,\,\,\, i=1,2.$$
 Assume $q_1<q<q_2$. Since
 $|\widetilde{u}|^q\leq
 |\widetilde{u}|^{q_1}+|\widetilde{u}|^{q_2}$, there holds
 $$\int_M |\widetilde{u}|^{q}dv_g\leq \int_M |\widetilde{u}|^{q_1}dv_g+\int_M |\widetilde{u}|^{q_2}dv_g
 \leq A_1^{q_1}+A_2^{q_2}.$$
 Hence
 $$\le(\int_M |{u}|^{q}dv_g\ri)^{1/q}\leq (A_1^{q_1}+A_2^{q_2})^{\f{1}{q}}\|u\|_{1,\tau}.$$
 Take $A=\max\{(A_1^{q_1}+A_2^{q_2})^{1/q_1},(A_1^{q_1}+A_2^{q_2})^{1/q_2}\}$. Then (\ref{concl}) follows
 immediately. $\hfill\Box$\\

 {\it Proof of proposition 2.1.} Assume there exist positive
constants $\alpha$, $\tau$ and $\beta$ such that (\ref{Tm-hold})
holds. For any $u\in W^{1,n}(M)$ we set
$\widetilde{u}=u/\|u\|_{1,\tau}$. It follows from (\ref{Tm-hold})
that
$$\int_M\sum_{k=n-1}^\infty\f{\alpha^k|\widetilde{u}|^{\f{nk}{n-1}}}{k!}dv_g\leq \beta.$$
Particularly for any integer $k\geq n-1$ there holds
$$\int_M\f{\alpha^k|\widetilde{u}|^{\f{nk}{n-1}}}{k!}dv_g\leq \beta,$$
and thus
$$\le(\int_M|u|^{\f{nk}{n-1}}dv_g\ri)^{\f{n-1}{nk}}\leq \le(\f{k!\beta}{\alpha^k}\ri)^{\f{n-1}{nk}}
\|u\|_{1,\tau}.$$
 For any $q\geq n$, there exists some $k\geq n-1$ such that
 $$\f{nk}{n-1}\leq q<\f{n(k+1)}{n-1}.$$
 In fact we can choose $k=[(n-1)p/n]$, the integer part of
 $(n-1)p/n$. By lemma 3.2, there exists a positive constant
 $A$ depending only on $n$, $q$, $\alpha$, and $\beta$ such that
 $$\le(\int_M |{u}|^{q}dv_g\ri)^{1/q}\leq A\|u\|_{1,\tau}.$$
 This implies that $W^{1,n}(M)\hookrightarrow
 L^q(M)$ continuously.
 Now we fix some $q>n$, say $q=n+1$. Then by lemma 3.1, there exists some constant $\epsilon>0$ depending only on
 $n$, $\alpha$, $\tau$, $\beta$ and $r$ such that for all $x\in M$,
 ${\rm Vol}_g(B_x(r))\geq \epsilon$.
 $\hfill\Box$\\

{\it Proof of corollary 2.2.} For any complete noncompact Riemannian
$n$-manifold
 $(M,g)$, if Trudinger-Moser inequality holds, then by proposition 2.1, there exists some constant $\epsilon>0$
 such that ${\rm Vol}_g(B_x(r))\geq \epsilon$ for all $x\in M$. Hence if there exists some
 complete noncompact
 Riemannian $n$-manifold $(M,g)$ such that
 $$\inf_{x\in M}{\rm Vol}_g(B_x(r))= 0,$$ then we conclude that Trudinger-Moser inequality is not valid on it.
 Now we construct such complete Riemannian manifolds. Consider the warped product
 $$M=\mathbb{R}\times N,\,\,\, g(t,\theta)=dt^2+f(t)ds_{N}^2,$$
 where $(N,ds_N^2)$ is a compact $(n-1)$-Riemannian manifold, $dt^2$ is the euclidian metric
 of $\mathbb{R}$, and $f$ is a smooth function satisfying
 $f(t)>0, \forall t\in\mathbb{R}$ and $\lim_{t\ra+\infty}f(t)= 0$.
 If $y=(t_1,m_1)$ and $z=(t_2,m_2)$ are two points of ${M}$,
 then $d_g(y,z)\geq |t_2-t_1|$. This together with the compactness of $N$ implies
 that $({M},g)$ is complete.
 In addition, for any $x=(t,m)\in M$, there holds
 $$B_x(1)\subset (t-1,t+1)\times N.$$
 Therefore
 \bea
  {\rm Vol}_g(B_{x}(1))&\leq&{\rm Vol}_g\le((t-1,t+1)\times
  N\ri){\nonumber}\\
  &\leq& {\rm
  Vol}_{ds_{N}^2}(N)\int_{t-1}^{t+1}f(t)dt{\nonumber}\\{\nonumber}
  &=&2{\rm Vol}_{ds_{N}^2}(N)f(\xi)\\{\label{volume}}&\ra& 0 \,\,\,{\rm as}\,\,\,
  t\ra+\infty,
 \eea
  where we used the integral mean value theorem, $\xi$ is some point in
  $(t-1,t+1)$. This gives the desired result.
   $\hfill\Box$

\section{Sufficient conditions}

In this section, we investigate sufficient conditions under which Trudinger-Moser inequality
holds. Precisely we shall prove theorem 2.3 and proposition 2.4. Firstly we have the following
key observation: \\

\noindent{\bf Lemma 4.1.} {\it Let $\mathbb{B}_0(\delta)\subset
\mathbb{R}^n$ be a ball centered at $0$ with radius $\delta$. If
$0\leq\alpha\leq \alpha_n=n\omega_{n-1}^{1/(n-1)}$, then there
exists some constant $C$ depending only on $n$ such that for all
$u\in W_0^{1,n}(\mathbb{B}_0(\delta))$ satisfying
$\int_{\mathbb{B}_0(\delta)}|\nabla u|^ndx\leq 1$, there holds
 \be\label{40}\int_{\mathbb{B}_0(\delta)}\zeta\le(n,\alpha|u|^{\f{n}{n-1}}\ri)dx
 \leq C\delta^n\le(\f{\alpha}{\alpha_n}\ri)^{n-1}
 \int_{\mathbb{B}_0(\delta)}|\nabla u|^ndx.\ee}

 \noindent{\it Proof.} Let $\widetilde{u}=u/\|\nabla
 u\|_{L^n(\mathbb{B}_0(\delta))}$. Since $\|\nabla
 u\|_{L^n(\mathbb{B}_0(\delta))}\leq 1$ and $0\leq\alpha\leq\alpha_n$, we have
 \bea
\zeta\le(n,\alpha|u|^{\f{n}{n-1}}\ri)&=&\sum_{k=n-1}^\infty
 \f{\alpha^k|u|^{\f{nk}{n-1}}}{k!}{\nonumber}\\
 &=&\sum_{k=n-1}^\infty\le(\f{\alpha}{\alpha_n}\ri)^{k}
 \f{\alpha_n^k\|\nabla
 u\|_{L^n(\mathbb{B}_0(\delta))}^{\f{nk}{n-1}}|\widetilde{u}|^{\f{nk}{n-1}}}{k!}{\nonumber}\\
 \label{4}
 &\leq&\|\nabla
 u\|_{L^n(\mathbb{B}_0(\delta))}^n\le(\f{\alpha}{\alpha_n}\ri)^{n-1}\zeta\le(n,\alpha_n|\widetilde{u}|^{\f{n}{n-1}}\ri).
 \eea
 It follows from the classical Trudinger-Moser inequality ((\ref{T-Mos}) with $\alpha_0$ replaced by $\alpha_n$)
  that
 \be\label{41}\int_{\mathbb{B}_0(\delta)}\zeta\le(n,\alpha_n|\widetilde{u}|^{\f{n}{n-1}}\ri)dx\leq
 C\delta^n\ee
 for some constant $C$ depending only on $n$. Integrating (\ref{4}) on $\mathbb{B}_0(\delta)$, we immediately
 obtain (\ref{40}) by using (\ref{41}). This concludes
 the lemma. $\hfill\Box$\\

 Let $(M,g)$ be a complete Riemannian $n$-manifold with
 ${\rm Ric}_{(M,g)}\geq Kg$ for some $K\in\mathbb{R}$ and ${\rm inj}_{(M,g)}\geq i_0$ for some
 $i_0>0$. Then we have the following local version of Trudinger-moser inequality which is the
 key estimate for the proof of theorem 2.3:\\

 \noindent{\bf Lemma 4.2.} {\it For any $\alpha: 0<\alpha<\alpha_n$ there exists some constant
 $\delta$
 depending only on $n$, $\alpha$, $K$ and $i_0$ such that for all $x\in M$
 and all $u\in C_0^{\infty}(B_x(\delta))$
 with $\|\nabla_g u\|_{L^n(B_x(\delta))}\leq 1$, there holds
 $$\int_M\zeta\le(n,\alpha|u|^{\f{n}{n-1}}\ri)dv_g\leq C\int_M|\nabla_g u|^ndv_g$$
 for some constant $C$ depending only on $n$, $\alpha$, $K$ and
 $i_0$.
 }\\

 \noindent{\it Proof.} By (Hebey \cite{Hebey}, theorem 1.3),
 we know that for any $\epsilon>0$
 there exists a positive constant $\delta$ depending only on $\epsilon$, $n$, $K$ and $i_0$ satisfying
 the following property: for any $x\in M$ there exists a harmonic
 coordinate chart $\phi: B_x(\delta)\ra\mathbb{R}^n$ such that
 $\phi(x)=0$, and the components $(g_{jl})$ of $g$ in this chart
 satisfy
 $$e^{-\epsilon}\delta_{jl}\leq g_{jl}\leq e^\epsilon \delta_{jl}$$
 as bilinear forms. One then has that $\phi(B_x(\delta))\subset
 \mathbb{B}_0(e^{\epsilon/2}\delta)$. Let $u$ be a function in
 $C_0^{\infty}(B_x(\delta))$ and $\|\nabla_g u\|_{L^n(B_x(\delta))}\leq 1$.
 It is not difficult to see that
 \bea\label{grad}\int_{B_x(\delta)}|\nabla_g u|^ndv_g&\geq&e^{-n\epsilon}\int_{\mathbb{B}_0(e^{\epsilon/2}\delta)}
 |\nabla(u\circ \phi^{-1})(x)|^ndx,\\\label{int}\int_M\zeta\le(n,\alpha|u|^{\f{n}{n-1}}\ri)dv_g
 &\leq&e^{n\epsilon/2}\int_{\mathbb{B}_0(e^{\epsilon/2}\delta)}
 \zeta\le(n,\alpha|(u\circ\phi^{-1})(x)|^{\f{n}{n-1}}\ri)dx.
 \eea
 For any fixed $\alpha: 0<\alpha<\alpha_n$, there exists some
 $\epsilon_0$ depending only on $n$ and $\alpha$ such that when $0<\epsilon\leq
 \epsilon_0$, it follows from (\ref{grad}) and $\|\nabla_g u\|_{L^n(B_x(\delta))}\leq
 1$ that
 $$\alpha\le(\int_{\mathbb{B}_0(e^{\epsilon/2}\delta)}|\nabla (u\circ\phi^{-1})(x)|^ndx\ri)^{1/(n-1)}
 \leq \alpha e^{n\epsilon_0/(n-1)}<\alpha_n.$$
  Now let $\epsilon=\epsilon_0$ be fixed and $\delta$ depending only on $\epsilon_0$, $n$, $K$ and $i_0$ be
 chosen as above. By lemma 4.1, there exists a constant $C_1=C_1(n)$ depending only on $n$
 such that
 $$\int_{\mathbb{B}_0(e^{\epsilon_0/2}\delta)}
 \zeta\le(n,\alpha|(u\circ\phi^{-1})(x)|^{\f{n}{n-1}}\ri)dx\leq C_1 e^{n\epsilon_0/2}\delta^n
 \int_{\mathbb{B}_0(e^{\epsilon_0/2}\delta)}|\nabla(u\circ\phi^{-1})(x)|^ndx.$$
 This together with (\ref{grad}) and (\ref{int}) implies that
 $$\int_M\zeta\le(n,\alpha|u|^{\f{n}{n-1}}\ri)dv_g\leq C_1e^{2n\epsilon_0}\delta^n\int_M|\nabla u|^ndv_g.$$
 Take $C=C_1e^{2n\epsilon_0}\delta^n$. We conclude that $C$ depends
 on $n$, $\alpha$, $K$ and $i_0$. $\hfill\Box$\\

 {\it Proof of theorem 2.3.} $(i)$ For any $\alpha: 0<\alpha<\alpha_n$, let
 $\delta=\delta(n,\alpha,K,i_0)$ be chosen as in lemma 4.2.
 Independently, by Gromov's covering lemma (Hebey \cite{Hebey}, lemma 1.6), we can select a sequence $(x_j)$ of
 points of $M$ such that\\

 \noindent $(a)$ $M=\cup_j B_{x_j}(\delta/2)$, and for any $ j\not=l$ there holds $B_{x_j}(\delta/4)\cap
 B_{x_l}(\delta/4)=\varnothing$;\\
 $(b)$ there exists $N$ depending only on $n$, $K$ and $\delta$ such that each
 point of $M$ has a neighborhood which intersects at most $N$ of the
 $B_{x_j}(\delta)$'s.\\

 \noindent For any $j$, we take a cut-off function $\phi_j\in
 C_0^\infty(B_{x_j}(\delta))$ satisfying $0\leq \phi_j\leq 1$, $\phi_j\equiv
 1$ on $B_{x_j}(\delta/2)$, and $|\nabla_g \phi_j|\leq 4/\delta$.
 It follows that for all $j$
 \be\label{9}|\nabla_g\phi_j^2|=2\phi_j|\nabla_g\phi_j|\leq \f{8}{\delta}\phi_j.\ee
 By the covering properties $(a)$ and $(b)$, we have
 \be\label{91}1\leq\sum_{j}\phi_j(x)\leq N\,\,\, {\rm for\,\,\,all}\,\,\, x\in M.\ee
 Set $\tau=8/\delta$. Assume $u\in C_0^\infty(M)$ satisfies
 $$\|u\|_{1,\tau}=\le(\int_M|\nabla u|^ndv_g\ri)^{1/n}+\tau
 \le(\int_M|u|^ndv_g\ri)^{1/n}\leq
 1.$$ It follows from (\ref{9}) and the Minkowvsky inequality that
 \bna\le(\int_M|\nabla_g(\phi_j^2 u)|^ndv_g\ri)^{1/n}\leq \le(\int_M\phi_j^{2n}|\nabla_g u|^ndv_g\ri)^{1/n}
 +\le(\int_M|\nabla_g\phi_j^2|^n|u|^ndv_g\ri)^{1/n}\leq\|u\|_{1,\tau}\leq 1.
 \ena
 In view of lemma 4.2, this leads to
 \bea
  \int_M\zeta\le(n,\alpha|u|^{\f{n}{n-1}}\ri)dv_g&\leq&\sum_{j}
 \int_{B_{\delta/2}(x_j)}\zeta\le(n,\alpha|u|^{\f{n}{n-1}}\ri)dv_g{\nonumber}\\{\nonumber}
 &\leq&\sum_{j}\int_{B_{\delta}(x_j)}\zeta\le(n,\alpha|\phi_j^2 u|^{\f{n}{n-1}}\ri)dv_g\\\label{16}
 &\leq&C\sum_{j}\int_M|\nabla(\phi_j^2 u)|^ndv_g
 \eea
 for some constant $C$ depending only on $n$, $\alpha$, $K$ and $i_0$.
 In addition we have by using
 (\ref{9}) and $0\leq \phi_j\leq 1$ that
 \bna
  \int_M|\nabla_g(\phi_j^2 u)|^ndv_g&\leq&
  2^n\int_M\le(\phi_j^{2n}|\nabla_g
  u|^n+|\nabla_g\phi_j^2|^n|u|^n\ri)dv_g{\nonumber}\\\label{26}
  &\leq&2^n\int_M\phi_j|\nabla_g
  u|^ndv_g+\f{16^n}{\delta^n}\int_M\phi_j|u|^ndv_g.
 \ena
 In view of (\ref{91}), it follows that
 \bna
  \sum_j\int_M|\nabla_g(\phi_j^2 u)|^ndv_g&\leq&2^nN\int_M|\nabla_g
  u|^ndv_g+\f{16^n}{\delta^n}N\int_M|u|^ndv_g\\
  &\leq&2^nN+\f{16^n}{\tau\delta^n}N.
 \ena
 This together with (\ref{16}) implies
 $$\int_M\zeta\le(n,\alpha|u|^{\f{n}{n-1}}\ri)dv_g\leq {C}$$
 for some constant ${C}$ depending only on $n$, $\alpha$, $K$ and $i_0$. By the density of $C_0^\infty(M)$ in $W^{1,n}(M)$,
 the inequality (\ref{Tm-hold}) holds for the above $\alpha$, $\tau$ and $C$.

 By proposition 2.1, we have that $W^{1,n}(M)$ is continuously embedded in $L^q(M)$ for any $q\geq n$.

 $(ii)$ Fix some point $z\in M$, let $r=r(x)=d_g(z,x)$ be the geodesic distance
 between $x$ and $z$. Without loss of generality, we may assume the injectivity radius of
 $(M,g)$ at $z$ is strictly larger than 1. Take a function sequence
 $$\phi_\epsilon(x)=\le\{\begin{array}{ll}
 1,&{\rm when}\quad r<\epsilon\\[1.5ex]
 \le(\log\f{1}{\epsilon}\ri)^{-1}\log\f{1}{r},&{\rm when}\quad \epsilon\leq
 r\leq 1\\[1.5ex] 0,&{\rm when}\quad r>1.
 \end{array}\ri.$$
 Then $\phi_\epsilon\in W^{1,n}(M)$ and for any constant $\tau>0$
 there holds
 \bna
 \|\phi_\epsilon\|_{1,\tau}=\le(\log\f{1}{\epsilon}\ri)^{(1-n)/n}\omega_{n-1}^{1/n}\le(1+O\le(\f{1}{\log{\epsilon}}\ri)\ri).
 \ena
 Set $\widetilde{\phi}_\epsilon=\phi_\epsilon/\|\phi_\epsilon\|_{1,\tau}$.
  Then we have on the geodesic ball $B_z(\epsilon)\subset M$,
 $$
 \zeta(n,\alpha\widetilde{\phi}_\epsilon^{\f{n}{n-1}})=
 e^{\alpha\widetilde{\phi}_\epsilon^{\f{n}{n-1}}}-\sum_{k=0}^{n-2}\f{\alpha^k\widetilde{\phi}_\epsilon^{\f{nk}{n-1}}}
 {k!}\geq
 \epsilon^{\alpha\omega_{n-1}^{-\f{1}{n-1}}(1+O(1/\log\epsilon))}+O\le(\le(\log\f{1}{\epsilon}\ri)^{n-2}\ri).
   $$
  Note that $\alpha\omega_{n-1}^{-\f{1}{n-1}}>n$ for any $\alpha>\alpha_n$. Hence, when $\alpha>\alpha_n$,
  we have
  \bna
   \int_M\zeta(n,\alpha|\widetilde{\phi}_\epsilon|^{\f{n}{n-1}})dv_g&\geq&\int_{B_z(\epsilon)}
   \zeta(n,\alpha|\widetilde{\phi}_\epsilon|^{\f{n}{n-1}})dv_g\\
   &\geq&\f{\omega_{n-1}}{n}(1+o_\epsilon(1))\epsilon^{n-\alpha\omega_{n-1}^{-1/(n-1)}(1+O(1/\log\epsilon))}+o_\epsilon(1).\\
   &\ra& +\infty\quad{\rm as}\quad \epsilon\ra 0.
  \ena
  This ends the proof of $(ii)$.

  $(iii)$  Take $\alpha_0: 0<\alpha_0<\alpha_n$. By $(i)$ there
  exists some $\tau_0=\tau_0(n,\alpha_0,K,i_0)>0$ such that
  $$\Lambda_{\alpha_0}:=\sup_{\|u\|_{1,\tau_0}\leq 1}\int_M
  \zeta(n,\alpha_0|u|^{\f{n}{n-1}})dv_g<\infty.$$
 Given any $\alpha>0$ and any $u\in W^{1,n}(M)$. Since $C_0^\infty(M)$ is dense in $W^{1,n}(M)$ under the norm
 $\|\cdot\|_{W^{1,n}(M)}$, which is equivalent to the norm
 $\|\cdot\|_{1,\tau_0}$, we can choose some $u_0\in C_0^\infty(M)$
 such that
 \be\label{alp0}2^{\f{n}{n-1}}\alpha\|u-u_0\|_{1,\tau_0}^{\f{n}{n-1}}<\alpha_0.\ee
 Since $\zeta(n,t)$ is increasing in $t$ for $t\geq 0$,  we obtain by using (\ref{Y2})
 \bea
  \int_M\zeta(n,\alpha|u|^{\f{n}{n-1}})dv_g&\leq&\int_M\zeta(n, 2^{\f{n}{n-1}}\alpha|u-u_0|^{\f{n}{n-1}}
  + 2^{\f{n}{n-1}}\alpha|u_0|^{\f{n}{n-1}})dv_g{\nonumber}\\{\nonumber}
  &\leq&\f{1}{\mu}\int_M\zeta(n,2^{\f{n}{n-1}}\alpha\mu|u-u_0|^{\f{n}{n-1}})dv_g\\\label{4-11}
  &&\quad+
  \f{1}{\nu}\int_M\zeta(n,2^{\f{n}{n-1}}\alpha\nu|u_0|^{\f{n}{n-1}})dv_g,
 \eea
 where $1/\mu+1/\nu=1$.  In view of (\ref{alp0}), we can take $\mu>1$
 sufficiently close to $1$
 such that
 $$2^{\f{n}{n-1}}\alpha\mu\|u-u_0\|_{1,\tau_0}^{\f{n}{n-1}}<\alpha_0.$$
 Hence
 \be\label{y2}\int_M\zeta(n,2^{\f{n}{n-1}}\alpha\mu|u-u_0|^{\f{n}{n-1}})dv_g\leq \Lambda_{\alpha_0}.\ee
 Since $u_0\in C_0^\infty(M)$, particularly $u_0$ has compact
 support, there holds
 \be\label{y3}\int_M\zeta(n,2^{\f{n}{n-1}}\alpha\nu|u_0|^{\f{n}{n-1}})dv_g<\infty.\ee
 Combining (\ref{4-11}), (\ref{y2}) and (\ref{y3}), we obtain
 $$\int_M\zeta(n,\alpha|u|^{\f{n}{n-1}})dv_g<\infty.$$
 This completes the proof of $(iii)$. $\hfill\Box$\\

 Now we shall prove proposition 2.4.
 Let us recall some notations from Riemannian geometry. In any chart,
 the Christoffel symbols of the Levi-Civita connection are given by
\be\label{connection}\Gamma_{ij}^k=\f{1}{2}g^{mk}\le(\p_ig_{mj}+\p_jg_{mi}-\p_mg_{ij}\ri),\ee
where $g_{ij}$'s are the components of $g$, $(g^{ij})$ denotes the
inverse matrix of $(g_{ij})$. Here and in the sequel the Einstein's
summation convention is adopted. Denote the Riemannian curvature of
$(M,g)$, a $(4,0)$-type tensor field, by ${\rm Rm}_{(M,g)}$. The
components of ${\rm Rm}_{(M,g)}$ are given by the relation
\be\label{curvature}R_{ijkl}=g_{i\alpha}\le(\p_k\Gamma_{jl}^\alpha-\p_l\Gamma_{jk}^\alpha
+\Gamma_{k\beta}^\alpha\Gamma_{jl}^\beta-\Gamma_{l\beta}^\alpha\Gamma_{jk}^\beta\ri).\ee
Similarly, the components of the Ricci curvature ${\rm Rc}_{(M,g)}$
of $(M,g)$ are given by the relation
\be\label{Ricci}R_{ij}=g^{\alpha\beta}R_{i\alpha j\beta}.\ee

 {\it Proof of proposition 2.4.} In view of proposition 2.1, it suffices to
 construct a complete noncompact Riemannian $n$-manifold $(M,g)$ such that its
 Ricci curvature has lower bound and there holds $$\inf_{x\in M}{\rm
 Vol}_g(B_x(1))=0.$$ Again we consider the warped product
 $${M}=\mathbb{R}\times N,\,\,\, g(x,\theta)=dx^2+f(x)ds_{N}^2,$$
 where $(N,ds_N^2)$ is a compact $(n-1)$-Riemannian manifold, $dx^2$ is the euclidean metric
 of $\mathbb{R}$, and $f$ is a smooth function satisfying
 $f(x)>0, \forall x\in\mathbb{R}$.
 In the following we calculate the Ricci curvature of
 $({M},g)$. In some product chart $(\mathbb{R}\times U,
 Id\times\phi)$ ($\{x,y^2,\cdots,y^n\}$), $g_{11}=1$,
 $g_{1\alpha}=0$, $g_{\alpha\beta}=fh_{\alpha\beta}$, $g^{11}=1$,
 $g^{1,\alpha}=0$, and
 $g^{\alpha\beta}=f^{-1}h^{\alpha\beta}$. Equivalently
 $$g=dx^2+f(x)h_{\alpha\beta}dy^\alpha dy^\beta,$$
 where $(h_{\alpha\beta})$ denote components of the metric $ds_N^2$.
  Here and in the sequel,  all indices $\alpha$, $\beta$,
 $\mu$, $\nu$ and $\lambda$ run from $2$ to $n$.
 In view of (\ref{connection}), the Christoffel symbols of the
 Levi-Civita connection was calculated as follows:
 \bna
 &&\Gamma_{11}^1=\Gamma_{11}^\alpha=\Gamma_{1\alpha}^1=0,\,\,\,
 \Gamma_{1\alpha}^\beta=\f{1}{2}g^{\mu\beta}\p_1g_{\mu\alpha}=\f{f^\prime}{2f}\delta_\alpha^\beta\\
 &&\Gamma_{\alpha\beta}^1=-\f{1}{2}\p_1g_{\alpha\beta}=-\f{f^\prime}{2}h_{\alpha\beta},\quad
 \Gamma_{\alpha\beta}^\gamma=\widetilde{\Gamma}_{\alpha\beta}^\gamma,
 \ena
 where $\delta_\alpha^\beta$ is equal to $1$ when $\alpha=\beta$, and $0$ when $\alpha\not=\beta$,
 $\widetilde{\Gamma}_{\alpha\beta}^\gamma$'s are
 components of the Christoffel symbols of
 Levi-Civita connection on $(N, ds_N^2)$.
 In view of (\ref{curvature}), the components of the Riemannian curvature reads
 \bna
  R_{1\alpha
  1\beta}&=&g_{11}\p_1\Gamma_{\alpha\beta}^1\\&=&\f{{f^\prime}^2-2ff^{\prime\prime}}{4f}h_{\alpha\beta}\\
  R_{1\alpha\beta\gamma}&=&g_{11}\le(\p_\beta\Gamma_{\alpha\gamma}^1-\p_\gamma\Gamma_{\alpha\beta}^1+
  \Gamma_{\beta k}^1\Gamma_{\alpha\gamma}^k-\Gamma_{\gamma
  k}^1\Gamma_{\alpha\beta}^k\ri)\\
  &=&\f{f^\prime}{2}\le(-\p_\beta h_{\alpha\gamma}+\p_\gamma h_{\alpha\beta}-h_{\beta\mu}\widetilde{\Gamma}_{\alpha\gamma}^\mu
  +h_{\gamma\mu}\widetilde{\Gamma}_{\alpha\beta}^\mu\ri)\\
  R_{\alpha\beta\gamma\mu}&=&g_{\alpha\lambda}\le(\p_\gamma\Gamma_{\beta\mu}^\lambda-\p_\mu\Gamma_{\beta\gamma}^\lambda
  +\Gamma_{\gamma k}^\lambda\Gamma_{\beta\mu}^k-\Gamma_{\mu
  k}^\lambda\Gamma_{\beta\gamma}^k\ri)\\
  &=&f\widetilde{R}_{\alpha\beta\gamma\mu}+g_{\alpha\lambda}\le(\Gamma_{\gamma 1}^\lambda\Gamma_{\beta\mu}^1-\Gamma_{\mu
  1}^\lambda\Gamma_{\beta\gamma}^1\ri)\\
  &=&f\widetilde{R}_{\alpha\beta\gamma\mu}+\f{{f^\prime}^2}{4}\le(h_{\alpha\mu}h_{\beta\gamma}-h_{\alpha\gamma}h_{\beta\mu}\ri),
 \ena
 where $\widetilde{R}_{\alpha\beta\gamma\mu}$'s denote the components
 of Riemannian curvature of $(N,ds_N^2)$.
 In view of (\ref{Ricci}), we get the components of the Ricci
 curvature as follows.
 \bna
  R_{11}&=&g^{\alpha\beta}R_{1\alpha
  1\beta}\\&=&(n-1)\f{{f^\prime}^2-2ff^{\prime\prime}}{4f^2}\\
  R_{1\alpha}&=&g^{\beta\gamma}R_{1\beta\alpha\gamma}\\
  &=&\f{f^\prime}{2f}h^{\beta\gamma}\le(-\p_\alpha h_{\beta\gamma}+\p_\gamma h_{\alpha\beta}
  -h_{\alpha\mu}\widetilde{\Gamma}_{\beta\gamma}^\mu+h_{\gamma\mu}\widetilde{\Gamma}_{\alpha\beta}^\mu\ri)\\
  R_{\alpha\beta}&=&g^{11}R_{\alpha1\beta1}+g^{\mu\nu}R_{\alpha\mu\beta\nu}\\
  &=&\f{{f^\prime}^2-2ff^{\prime\prime}}{4f}h_{\alpha\beta}+\widetilde{R}_{\alpha\beta}+
  \f{{f^\prime}^2}{4f}h^{\mu\nu}\le(h_{\alpha\nu}h_{\mu\beta}-h_{\alpha\beta}h_{\mu\nu}\ri)\\
  &=&\f{(2-n){f^\prime}^2-2ff^{\prime\prime}}{4f}h_{\alpha\beta}+\widetilde{R}_{\alpha\beta},
 \ena
 where $\widetilde{R}_{\alpha\beta}$'s are components of the Ricci
 curvature of $(N,ds_N^2)$.
 If we assume the functions $f$, $f^\prime/f$
 and $f^{\prime\prime}/f$ are all bounded, then in the chart $(\mathbb{R}\times U,
 Id\times\phi)$, the eigenvalues of the matrix $(R_{jl})$ and the matrix $(g_{jl})$ are uniformly bounded.
 Thus there exists some constant $K_1\in\mathbb{R}$ such that $(R_{jl})\geq K_1(g_{jl})$.
   Note that $(N,ds_N^2)$ is compact. There exists some constant $K\in\mathbb{R}$ such that
    $Ric_{(M,g)}\geq Kg$
 as bilinear forms. If we further assume
 $\lim_{x\ra+\infty}f(x)=0$, then by (\ref{volume}), we have ${\rm Vol}_g(B_y(1))\ra
 0$ as $x\ra+\infty$, where $y=(x,m)\in \mathbb{R}\times N$. One can check that the following functions satisfy
 all the above assumptions on $f$.

 \noindent $\bullet$ $f$ is a smooth positive function defined on
 $\mathbb{R}$ and satisfies
 $$f(x)=\le\{\begin{array}{lll}
  (1+x^2)e^{-x+\sin x}\quad&{\rm when}\quad &x>1\\[1.5ex]
  1, &{\rm when}& x<0
 \end{array}\ri.$$
 $\bullet$ $f$ is a smooth positive function defined on
 $\mathbb{R}$ and satisfies
 $$f(x)=\le\{\begin{array}{lll}
  \f{1}{\log x}\quad&{\rm when}\quad &x>2\\[1.5ex]
  1, &{\rm when}& x<0.
 \end{array}\ri.$$
  This gives the desired result.$\hfill\Box$

\section{Proof of proposition 2.5}
In this section, we shall construct complete noncompact Riemannain
$n$-manifolds to show that the condition Ricci curvature has lower
bound in theorem 2.3 is not necessarily needed.\\

{\it Proof of proposition 2.5.} It suffices to construct a complete
noncompact Riemannian $n$-manifold on which Trudinger-Moser
embedding holds, but its Ricci curvature has no lower bound. For
this purpose, we consider the Riemannian manifold
$(\mathbb{R}^n,g)$, where $\mathbb{R}^n$ is the euclidian space and
$$g=dx_1^2+f(x_1)dx_2^2+\cdots+f(x_1)dx_n^2,$$ and $f$ is a smooth function
on $\mathbb{R}$ such that $a\leq f\leq b$ for two positive constants
$a$ and $b$. Clearly $(\mathbb{R}^n,g)$ is complete and noncompact.
In view of Trudinger-Moser inequality on the standard euclidian
space $\mathbb{R}^n$ \cite{Cao,doo1,Panda}, one can easily see that
if $\alpha$ is chosen sufficiently small, then the supremum
$$\sup_{u\in W^{1,n}(M),\,\|u\|_{W^{1,n}}\leq 1}\int_{\mathbb{R}^n}\zeta\le(n,\alpha|u|^{\f{n}{n-1}}\ri)dv_g$$
 is finite, i.e. Trudinger-Moser inequality holds on the manifold
 $(\mathbb{R}^n,g)$, where
 $$\|u\|_{W^{1,n}}=\le(\int_{\mathbb{R}^n}(|\nabla_gu|^n+|u|^n)dv_g\ri).$$

In the following, we shall further choose $f$ such that the Ricci
curvature of $(\mathbb{R}^n,g)$ is unbounded from below. By
(\ref{Ricci}),
\be\label{rc}R_{11}=(n-1)\f{{f^\prime}^2-2ff^{\prime\prime}}{4f^2}.\ee
It suffices to find a sequence of points $(x^{(m)})$ of
$\mathbb{R}^n$ such that $R_{11}(x^{(m)})\ra -\infty$. One choice of
$f$ is that $f(t)=2+\sin t^2$. In this case, we have \bna
f^\prime(x_1)=2+2x_1 \cos x_1^2,\quad f^{\prime\prime}(x_1)=2\cos
x_1^2-4x_1^2\sin x_1^2.\ena Thus (\ref{rc}) implies
$$R_{11}(x)=(n-1)\f{(2+2x_1\cos x_1^2)^2-2(2+\sin x_1^2)(2\cos x_1^2-4x_1^2\sin x_1^2)}
{4(2+\sin x_1^2)^2}.$$ Choosing
$x^{(m)}=\le(\sqrt{2m\pi+3\pi/2},0,\cdots,0\ri)$, we obtain
$$R_{11}({x^{(m)}})=-4m\pi-3\pi+n-1\ra -\infty\quad{\rm as}\quad
m\ra\infty.$$ Another choice of $f$ is that $f(t)=e^{\sin t^2}$. In
this case, we have
\bna
 f^\prime(x_1)&=&2x_1e^{\sin x_1^2}\cos x_1^2,\quad
 f^{\prime\prime}(x_1)=e^{\sin x_1^2}\le(-4x_1^2\sin x_1^2+4x_1^2\cos^2x_1^2+2\cos
 x_1^2\ri).
\ena In view of (\ref{rc}), we obtain
$$R_{11}(x)=(n-1)(2x_1^2\sin x_1^2+x_1^2\cos^2x_1^2-2x_1^2\cos^2 x_1^2-\cos x_1^2).$$
Again, we select $x^{(m)}=\le(\sqrt{2m\pi+3\pi/2},0,\cdots,0\ri)$
and conclude $R_{11}(x^{(m)})\ra -\infty$ as $m\ra\infty$.
$\hfill\Box$

\section{Adams inequalities}
In this section, we concern Adams inequalities on complete noncompact Riemannian manifolds.
Precisely we shall prove theorem 2.6. The method we adopted
here is similar to that of theorem 2.3. \\

 \noindent{\it Proof of theorem 2.6.} $(i)$ Suppose that ${\rm
 inj}_{(M,g)}\geq i_0>0$ and there exist constants $C(k)$ such that $|\nabla^k{\rm Rc}_{(M,g)}|\leq
 C(k)$, $k=0, 1,\cdots, m-1$. It follows from (Hebey \cite{Hebey}, theorem 1.3)
 that for any $Q>1$ and
 $\alpha\in(0,1)$, the harmonic radius $r_H=r_H(Q,m,\alpha)$ is
 positive. Namely, for any $Q>1$, $\alpha\in(0,1)$, and $x\in M$, there
 exists a harmonic coordinate chart $\psi: B_x(r_H)\ra \mathbb{R}^n$ such
 that
 \be\label{cmcontroll}
 \le\{
 \begin{array}{lll}
  Q^{-1}\delta_{lq}\leq g_{lq}\leq Q\delta_{lq}\quad {\rm as\,\, a\,\,
 bilinear\,\,
 form};\\[1.5ex]
 \sum_{1\leq|\beta|\leq m}\|\p^\beta g_{lq}\|_{C^0(B_x(r_H))}+
 \sum_{|\beta|=m}\|\p^\beta g_{lq}\|_{C^\alpha(B_x(r_H))}\leq
 Q-1.\end{array}\ri.\ee
 Now we fix $Q>1$ and $\alpha\in(0,1)$. Without loss of generality, we may assume $\psi(x)=0$.
 Particularly we have that for any $r:0<r\leq r_H$
 $$\mathbb{B}_0(r/\sqrt{Q})\subset\psi(B_x(r))\subset\mathbb{B}_0({\sqrt{Q}}r).$$
 Let $\eta\in C_0^\infty(\mathbb{R}^n)$ be such that
 $0\leq\eta\leq 1$, and
  $$\eta=\le\{\begin{array}{lll}1&\,\,\,{\rm on}\,\,\, &\mathbb{B}_0(r_H/(4\sqrt{Q})),\\[1.5ex]
  0&\,\,\,{\rm on}&\mathbb{R}^n\setminus\mathbb{B}_0(r_H/(2\sqrt{Q})).\end{array}\ri.$$
  Then $\eta\circ\psi\in
 C_0^\infty(M)$ satisfies $0\leq \eta\circ\psi\leq 1$, $\eta\circ\psi\equiv
 1$ on $B_x(r_H/(4Q))$, and $\eta\circ\psi\equiv 0$ on $M\setminus
 B_x(r_H/2)$. By Gromov's covering lemma (Hebey \cite{Hebey}, lemma 1.6), there exists
 a sequence of points $(x_k)$  of $M$
 such that
 \be\label{cover}M=\cup_{k}B_{x_k}(r_H/(4Q))\ee
 and there exists some integer $N$ such that for any $x\in M$, $x$ belongs to
 at most $N$ balls in the covering.
 Let $\psi_k: B_{x_k}(r_H)\ra \mathbb{R}^n$ be as the above $\psi$ and set
 $\eta_k=\eta\circ\psi_k$. By (\ref{cmcontroll}), the
 components of the metric tensor are $C^{m}$-controlled in the
 charts $(B_{x_k}(r_H),\psi_k)$. It then follows that there exists some constant $C_1>0$
 depending only on $r_H$ and $Q$
 such that $|\nabla_g^l\eta_k|\leq
 C_1$ for all all $l:0\leq l\leq m$ and all $k\in \mathbb{N}$, where $\nabla_g^l$ is defined by (\ref{gradient}).

 Assume $u\in C^\infty(M)$ satisfies $\|u\|_{W^{m,n/m}(M)}\leq 1$.
 Then we get $$\eta_k^{m+1}u\in C_0^\infty(B_{x_k}(r_H/2))$$ and
 \be\label{etak}\|\nabla_g^m(\eta_k^{m+1}u)\|_{L^{\f{n}{m}}(B_{x_k}(r_H/2))}\leq C_2\ee
 for some constant $C_2$ depending only on $n$, $m$, and $C_1$.
 By the standard elliptic estimates (Gilbarg-Trudinger
 \cite{GT}, Chapter 9), one can see that
 \be\label{gra}\|\nabla_{\mathbb{R}^n}^m((\eta_k^{m+1}u)\circ\psi_k^{-1})\|_{L^{\f{n}{m}}(\mathbb{B}_0(\sqrt{Q}r_H))}\leq C_3\ee
 for some constant $C_3$ depending only on $n$, $m$, $Q$, $r_H$ and $C_1$.
 Let $j$ be the smallest
 integer great than or equal to $n/m$. Similarly as we derived (\ref{16}), we
 calculate by using (\ref{cover}), (\ref{etak}) and the relation ${(j-1)n}/{(n-m)}\geq {n}/{m}$
 \bea
 \int_M\zeta\le(j,\alpha|u|^{\f{n}{n-m}}\ri)dv_g&\leq&\sum_k\int_{B_{x_k}(r_H/(4Q))}
 \zeta\le(j,\alpha|u|^{\f{n}{n-m}}\ri)dv_g{\nonumber}\\
 &\leq&\sum_k\int_{B_{x_k}(r_H/2)}\zeta\le(j,\alpha|\eta_k^{m+1}u|^{\f{n}{n-m}}\ri)dv_g{\nonumber}\\
 &\leq&\sum_k\le(\f{\|\nabla_g^m(\eta_k^{m+1}
 u)\|_{L^{\f{n}{m}}(B_{x_k}(r_H/2))}}{C_2}\ri)^{\f{(j-1)n}{n-m}}\int_{B_{x_k}(r_H/2)}
 \zeta\le(j,\alpha
 C_2^{\f{n}{n-m}}|\eta_k^{m+1}u|^{\f{n}{n-m}}\ri)dv_g{\nonumber}\\\label{mf}
 &\leq&\sum_k\f{\|\nabla_g^m(\eta_k^{m+1}
 u)\|_{L^{\f{n}{m}}(B_{x_k}(r_H/2))}^{\f{n}{m}}}{C_2^{\f{n}{m}}}\int_{B_{x_k}(r_H/2)}
 \zeta\le(j,\alpha C_2^{\f{n}{n-m}}|\eta_k^{m+1}u|^{\f{n}{n-m}}\ri)dv_g.
 \eea
  Noting that $Q^{-1}\delta_{lq}\leq g_{lq}\leq Q\delta_{lq}$ as a bilinear form, we
  have
 \be\label{rh}\int_{B_{x_k}(r_H/2)}
 \zeta\le(j,\alpha C_2^{\f{n}{n-m}}|\eta_k^{m+1}u|^{\f{n}{n-m}}\ri)dv_g\leq Q^{\f{n}{2}}
 \int_{\mathbb{B}_0(\sqrt{Q}r_H)}
 \zeta\le(j,\alpha
 C_2^{\f{n}{n-m}}|(\eta_k^{m+1}u)\circ\psi_k^{-1}|^{\f{n}{n-m}}\ri)dx.\ee
 In view of (\ref{gra}), we take
 \be\label{alpha-0}\alpha_0=\beta_0/(C_2C_3)^{\f{n}{n-m}}.\ee
 Then for any
 $\alpha: 0<\alpha\leq\alpha_0$, it follows from Adams inequality (\ref{Ada}) that
 \be\label{e-u}\int_{\mathbb{B}_0(\sqrt{Q}r_H)}
 \zeta\le(j,\alpha C_2^{\f{n}{n-m}}|(\eta_k^{m+1}u)\circ\psi_k^{-1}|^{\f{n}{n-m}}\ri)dx
 \leq C_{m,n}|\mathbb{B}_0(\sqrt{Q}r_H)|.\ee
 Clearly there exists some constant $C_4>0$ depending only on $n$, $m$, $Q$ and $r_H$ such that
 \be\label{gra-eta}|\nabla_g^l\eta_k^{m+1}|^{\f{n}{m}}\leq C_4\eta_k,\quad\forall l=0,1,\cdots,m.\ee
 Since $1\leq\sum_k\eta_k(x)\leq N$ for all $x\in M$, we obtain by combining
 (\ref{mf})-(\ref{gra-eta}) that
 \bna\int_M\zeta\le(j,\alpha|u|^{\f{n}{n-m}}\ri)dv_g&\leq&
 C_5\sum_k\int_M|\nabla_g^m(\eta_k^{m+1}
 u)|^{\f{n}{m}}dv_g\\
 &\leq& C_5\sum_k\sum_{l=0}^m(C_m^l)^{\f{n}{m}}\int_M|\nabla_g^{m-k}\eta_k^{m+1}\nabla_g^lu|^{\f{n}{m}}dv_g\\
 &\leq&C_4C_5\sum_{l=0}^m(C_m^l)^{\f{n}{m}}\int_M(\sum_k\eta_k)|\nabla_g^lu|^{\f{n}{m}}dv_g\\
 &\leq&C_4C_5N\sum_{l=0}^m(C_m^l)^{\f{n}{m}}\int_M|\nabla_g^lu|^{\f{n}{m}}dv_g\\
 &\leq& C_6\ena
 for constants $C_5$ and $C_6$ depending only on $n$, $m$, $Q$ and
 $r_H$, where $C_m^l=\f{m!}{l!\,(m-l)!}$.

 According to (Hebey \cite{Hebey}, theorem 2.8),
 $C_0^\infty(M)$ is dense in $W^{m,\f{n}{m}}(M)$. Hence
for any $u\in W^{m,\f{n}{m}}(M)$, there exists a sequence $(u_k)$ in
$C_0^\infty(M)$ such that $\|u_k-
u\|_{W^{m,\f{n}{m}}(M)}\ra 0$ as $k\ra\infty$. Assume
$\|u\|_{W^{m,\f{n}{m}}(M)}\leq 1$. Then for any $\alpha:
0<\alpha<\alpha_0$ there holds \bna
 \int_M\zeta\le(j,\alpha|u|^{\f{n}{n-m}}\ri)dv_g&\leq&\lim_{k\ra\infty}\int_M\zeta
 \le(j,\alpha|u_k|^{\f{n}{n-m}}\ri)dv_g\leq
 C_6.
\ena
 Using the same method of deriving $W^{1,n}(M)\hookrightarrow
 L^q(M)$ continuously for all $q\geq n$ in theorem 2.3, we obtain
 the continuous embedding
 $W^{m,{n}/{m}}(M)\hookrightarrow L^q(M)$ for any $q\geq {n}/{m}$.\\

 $(ii)$ Let $\alpha>0$ be any real number and
 $u$ be any function belonging to the space
 $W^{m,\f{n}{m}}(M)$. Since $C_0^\infty(M)$ is dense in $W^{m,\f{n}{m}}(M)$,
  there exists
 some $u_0\in C_0^\infty(M)$ such that
 \be\label{cd}\alpha\|u-u_0\|_{W^{m,\f{n}{m}}(M)}^{\f{n}{n-m}}<\alpha_0/2,\ee
 where $\alpha_0$ is defined by (\ref{alpha-0}).
 Using (\ref{Y2}) and an elementary inequality
 $$|a|^p\leq (1+\epsilon)|a-b|^p+c(\epsilon,p)|b|^p,$$ where
 $\epsilon>0$, $p>1$ and $c(\epsilon,p)$ is a constant depending only on $\epsilon$ and
 $p$, we have
 \bea
  \int_M\zeta\le(j,\alpha|u|^{\f{n}{n-m}}\ri)dv_g&\leq&\int_M\zeta\le(j,(1+\epsilon)\alpha|u-u_0|^{\f{n}{n-m}}+
  c(\epsilon,n/(n-m))\alpha|u_0|^{\f{n}{n-m}}\ri)dv_g\nonumber\\\nonumber
  &\leq&\f{1}{\mu}\int_M\zeta\le(j,\mu
  (1+\epsilon)\alpha|u-u_0|^{\f{n}{n-m}}\ri)dv_g\\\label{u0}
  &&+\f{1}{\nu}\int_M\zeta\le(j,\nu
  c(\epsilon,n/(n-m))\alpha|u_0|^{\f{n}{n-m}}\ri)dv_g,
 \eea
 where $\mu>1$, $\nu>1$ and $1/\mu+1/\nu=1$. Choosing $\epsilon$
 sufficiently small and $\mu$ sufficiently close to $1$ such that
 $\mu(1+\epsilon)\alpha_0/2\leq \alpha_0$, in view of (\ref{cd}), we
 have by part $(i)$
 \be\label{6-11}\int_M\zeta\le(j,\mu
  (1+\epsilon)\alpha|u-u_0|^{\f{n}{n-m}}\ri)dv_g\leq C_6.\ee
   Note that $u_0\in
  C_0^\infty(M)$, particularly $u_0$ has compact support. It follows
  that
  \be\label{6-12}\int_M\zeta\le(j,\nu
  c(\epsilon,n/(n-m))\alpha|u_0|^{\f{n}{n-m}}\ri)dv_g<\infty.\ee
  Inserting (\ref{6-11}) and (\ref{6-12}) into (\ref{u0}), we complete the proof of part $(ii)$. $\hfill\Box$

\section{Applications of Trudinger-Moser inequalities}

 In this section, we consider applications of theorem 2.3, namely the existence and multiplicity results for the problem
 (\ref{equa}) and its perturbation (\ref{equa1}). Specifically we shall prove theorem 2.7 and theorem 2.10.
 Throughout this section, we use the notations introduced in section 2. Let $(M,g)$ be a complete noncompact
 Riemannian $n$-manifold with ${\rm Rc}_{(M,g)}\geq Kg$ for some $K\in\mathbb{R}$ and ${\rm inj}_{(M,g)}\geq
 i_0>0$. Assume $\phi(x)$ satisfies the hypotheses $(\phi_1)$ and $(\phi_2)$,
 $v(x)$ satisfies the hypotheses $(v_1)$ and $(v_2)$. Let $E$ be a function space defined by
 (\ref{space}). If $u\in E$, then the $E$-norm of $u$ is defined by
 $$\|u\|_E=\le(\int_M(|\nabla_g u|^n+v|u|^n)dv_g\ri)^{1/n}.$$

  The following compact embedding result is very important in our analysis.\\

 \noindent{\bf Proposition 7.1.} {\it For any $q\geq n$, the function space $E$ is compactly embedded in $L^q(M)$.}\\

 \noindent{\it Proof.} Let $(u_k)$ be a sequence of functions with $\|u_k\|_{E}\leq C$ for some constant
 $C$. It suffices to prove that up to a subsequence, $(u_k)$ converges in $L^q(M)$ for any $q\geq n$.
 Clearly $(u_k)$ is bounded in $W^{1,n}(M)$, and thus we can
 assume that for any $q>1$, up to a subsequence
 \bea
  &&u_k\rightharpoonup u_0\quad{\rm weakly\,\,in}\quad
  E{\nonumber}\\{\label{Lploc}}
  &&u_k\ra u_0\quad{\rm strongly\,\,in}\quad L^q_{\rm loc}(M)\\
  &&u_k\ra u_0\quad{\rm a.\,e.}\,\,\, {\rm in}\quad M.\nonumber
 \eea
  If $v(x)\in L^{1/(n-1)}(M)$, using
 the same argument of (\cite{Yang}, Lemma 2.4), we conclude that $E\hookrightarrow L^q(M)$
 compactly for any $q>1$.
 So, in view of $(v_2)$, we may assume $v(x)\ra\infty$ as $d_g(O,x)\ra\infty$,
 where $O$ is a fixed point of $M$.
  Given any $\epsilon>0$, there
 exists some $R>0$ such that $v(x)>(2C)^n/\epsilon$ when
 $d_g(O,x)\geq R$. Hence
 $$\f{(2C)^n}{\epsilon}\int_{M\setminus B_O(R)}|u_k-u_0|^ndv_g<\int_Mv|u_k-u_0|^ndv_g\leq (2C)^n.$$
 This gives
 $$\int_{M\setminus B_O(R)}|u_k-u_0|^ndv_g<\epsilon.$$
 By (\ref{Lploc}), we have
 $$\lim_{k\ra\infty}\int_{B_O(R)}|u_k-u_0|^ndv_g=0.$$
 Hence for the above $\epsilon$, there exists some
 $l\in\mathbb{N}$ such that when $k>l$,
 $$\int_M|u_k-u_0|^ndv_g<2\epsilon.$$
 This implies $u_k\ra u_0$ strongly in $L^n(M)$ as $k\ra\infty$.

 It follows from $(i)$ of theorem 2.3 that $(u_k)$ is bounded in $L^q(M)$ for any $q\geq n$.
 Now fixing $q>n$, we get by H\"older's inequality
 $$\int_M|u_k-u_0|^qdv_g\leq\le(\int_M|u_k-u_0|^ndv_g\ri)^{1/n}
 \le(\int_{M}|u_k-u_0|^{\f{n(q-1)}{n-1}}dv_g\ri)^{1-1/n}.$$
 This together with the fact that $u_k\ra u_0$ in $L^n(M)$
 implies $u_k\ra u_0$ in $L^q(M)$.
 $\hfill\Box$\\

 Let $S_q$ be defined by (\ref{Sq}). Then we have the following:\\

 \noindent{\bf Proposition 7.2.} {\it For any $q>n$, $S_q$ is attained by some nonnegative function $u\in
 E\setminus\{0\}$.}\\

 \noindent{\it Proof.} Assume $q>n$. It is easy to see that
 $$S_q^n=\inf_{\int_M\phi|u|^qdv_g=1}\int_{M}\le(|\nabla u|^n+v|u|^n\ri)dv_g.$$
 Choosing a sequence of functions $(u_k)\subset
 E$ such that $\int_M\phi|u_k|^qdv_g=1$ and
 $$\lim_{k\ra\infty}\int_{M}\le(|\nabla u_k|^n+v|u_k|^n\ri)dv_g=S_q^n.$$
 By proposition 7.1, there exists some $u\in E$ such that up to a subsequence, $u_k\rightharpoonup u$
 weakly in $E$, $u_k\ra u$ strongly in $L^q(M)$ for any $q\geq
 n$, and $u_k\ra u$ almost everywhere in $M$.
 Since $u_k\ra u$ strongly in $L^s(B_O(R_0))$ for all $s>1$ and $\phi\in L^p(B_O(R_0))$, we have
 by using H\"older's inequality that
 \be\label{ho1}\lim_{k\ra \infty}\int_{B_O(R_0)}\phi|u_k|^qdv_g=\int_{B_O(R_0)}\phi|u|^qdv_g.\ee
 In view of $(v_2)$, we have
 \bna
  \int_{M\setminus B_O(R_0)}\phi||u_k|^q-|u|^q|dv_g&\leq&
  qC_0\int_M(|u_k|^{q-1}+|u|^{q-1})|u_k-u|dv_g\\
  &\leq&
  qC_0\le\{\le(\int_M|u_k|^qdv_g\ri)^{1-1/q}+\le(\int_M|u|^qdv_g\ri)^{1-1/q}\ri\}\\
  &&\quad \times\le(\int_M|u_k-u|^qdv_g\ri)^{1/q}\\
  &\ra& 0\,\,\,{\rm as}\,\,\, k\ra\infty.
 \ena
 This together with (\ref{ho1}) implies
 \be\label{L1}\int_M\phi|u|^qdv_g=\lim_{k\ra\infty}\int_M\phi|u_k|^qdv_g=1.\ee
 Since $u_k\rightharpoonup u$ weakly in $E$, we have
 \bna
 \int_M|\nabla u|^ndv_g=\lim_{k\ra\infty}\int_M|\nabla
 u|^{n-2}\nabla u\nabla u_kdv_g\leq\limsup_{k\ra\infty}\le(\int_M|\nabla u_k|^ndv_g\ri)^{\f{1}{n}}
 \le(\int_M|\nabla u|^ndv_g\ri)^{1-\f{1}{n}},
 \ena
 from which we obtain
 \be\label{L2}\int_M|\nabla u|^ndv_g\leq \limsup_{k\ra\infty}\int_M|\nabla
 u_k|^ndv_g.\ee
 In addition, we have by Fatou's lemma
 \be\label{L3}\int_Mv|u|^ndv_g\leq
 \limsup_{k\ra\infty}\int_Mv|u_k|^ndv_g.\ee
 Combining (\ref{L1}), (\ref{L2}) and (\ref{L3}), we conclude that $S_q$ is attained by $u\in E\setminus
 \{0\}$. Since $|u|\in E$, one can easily see that $S_q$ is also attained by $|u|$. $\hfill\Box$\\

  Now we get back to the problem (\ref{equa}). Since we are interested in nonnegative weak solutions,
 without loss of generality we may
 assume $f(x,t)\equiv 0$ for all $(x,t)\in M\times(-\infty,0]$. By $(f_1)$, we have for all
 $(x,t)\in M\times\mathbb{R}$,
 $$|F(x,t)|\leq \f{b_1}{n}|t|^n+b_2t\zeta\le(n,|t|^{\f{n}{n-1}}\ri).$$
 This together with $(\phi_1)$, $(\phi_2)$ and (\ref{Y1}) implies
 that for any $u\in E$ there holds
 \bna
  \int_M\phi
  F(x,u)dv_g&\leq&\|\phi\|_{L^p(B_O(R_0))}\|F(x,u)\|_{L^q(M)}+C_0\int_MF(x,u)dv_g\\
  &\leq&\|\phi\|_{L^p(B_O(R_0))}\le(\f{b_1}{n}\|u\|_{L^{qn}(M)}^n+b_2\|u\zeta(n,|u|^{\f{n}{n-1}})\|
  _{L^q(M)}\ri)\\
  &&+C_0\f{b_1}{n}\|u\|_{L^n(M)}^n+C_0b_2\|u\zeta(n,|u|^{\f{n}{n-1}})\|_{L^1(M)}\\
  &\leq&C\le(\|u\|_{L^{qn}(M)}^n+\|u\|_{L^{qn}(M)}\|\zeta(n,\f{qn}{n-1}|u|^{\f{n}{n-1}})\|_{L^1(M)}
  ^{1-\f{1}{n}}\ri.\\
  &&\le.
  \|u\|_{L^n(M)}^n+\|u\|_{L^n(M)}\|\zeta(n,\f{n}{n-1}|u|^{\f{n}{n-1}})\|_{L^1(M)}\ri),
 \ena
 where $C$ is a constant depending only on $n$, $b_1$, $b_2$, $C_0$
 and  $\|\phi\|_{L^p(B_O(R_0))}$, and $1/p+1/q=1$. By
 theorem 2.3, $u\in L^s(M)$ for all $s\geq n$, and for any $\alpha>0$ there holds $\zeta(n,\alpha|u|^{\f{n}{n-1}})\in
 L^1(M)$.  Hence
 $$\int_M\phi
  F(x,u)dv_g<+\infty,\quad\forall u\in E.$$
  Based on this, we can define a functional on $E$ by
 \be\label{J}J(u)=\f{1}{n}\|u\|_E^n-\int_M\phi F(x,u)dv_g.\ee
 By (\cite{doo}, proposition 1) and the standard argument \cite{Rabin},
 we have $J\in \mathcal{C}^1(E,\mathbb{R})$. Clearly the critical point of $J$ is a weak solution to
 (\ref{equa}). Concerning the geometry of $J$, the following two
 lemmas imply that $J$ has a mountain pass structure.\\

 \noindent{\bf Lemma 7.3.} {\it Assume that $(f_1)$, $(f_2)$, and $(f_3)$ are satisfied.
  Then for any nonnegative, compactly supported function $u\in
  E\setminus\{0\}$, there holds $J(tu)\ra
  -\infty$ as $t\ra +\infty$.}\\

  \noindent{\it Proof.} By $(f_2)$ and $(f_3)$, there exist $c_1$, $c_2>0$ and $\mu>n$ such that $F(x,s)\geq
  c_1s^\mu-c_2$ for all $(x,s)\in M\times [0,+\infty)$. Assume ${\rm supp}\,u\subset B_O(R_1)$ for some
  $R_1>0$.
  We have
  \bna
   J(t u)&=&\f{t^n}{n}\|u\|_E^n-\int_{B_O(R_1)}\phi F(x,t
   u)dv_g\\&\leq&\f{t^n}{n}\|u\|_E^n-c_1
   t^\mu\int_{B_O(R_1)}\phi u^\mu dv_g-c_2\int_{B_O(R_1)}\phi dv_g.
  \ena
   This gives the desired result since $\phi(x)>0$ for all $x\in M$ and $\mu>n$.$\hfill\Box$\\

 \noindent {\bf Lemma 7.4.} {\it Assume that
 $(f_1)$ and $(f_4)$ are satisfied. Then there exist sufficiently small constants $r>0$ and $\delta>0$
  such that $J(u)\geq\delta$ for all $u$ with $\|u\|_E=r$.}\\

  \noindent {\it Proof.}
  By $(f_1)$ and $(f_4)$, there exists some constants $\theta\in (0,1)$
  and $C>0$ such that
  $$F(x,s)\leq \f{(1-\theta)\lambda_\phi}{n}|s|^n+C|s|^{n+1}\zeta\le(n,\alpha_0|s|^{\f{n}{n-1}}\ri)$$
  for all $(x,s)\in M\times\mathbb{R}$. By definition of
  $\lambda_\phi$,
  \be\label{lam-p}\f{(1-\theta)\lambda_\phi}{n}\int_M\phi|u|^ndv_g\leq
  \f{1-\theta}{n}\|u\|_E^n.\ee
  Note that $\phi$ satisfies $(\phi_1)$ and $(\phi_2)$. We have by H\"older's inequality
  and (\ref{Y1}) that
  \bea\nonumber
  \int_M\phi|u|^{n+1}\zeta\le(n,\alpha_0|u|^{\f{n}{n-1}}\ri)dv_g&\leq&\|\phi\|_{L^p(B_O(R_0))}
  \le(\int_M|u|^{(n+1)q}dv_g\ri)^{1/q}\le(\int_M\zeta\le(n,q^\prime\alpha_0|u|^{\f{n}{n-1}}\ri)dv_g\ri)^{1/q^\prime}\\
  \label{mmm}&&+C_0\le(\int_M|u|^{(n+1)\beta}dv_g\ri)^{1/\beta}\le(\int_M\zeta\le(n,\gamma\alpha_0|u|^{\f{n}{n-1}}\ri)dv_g\ri)^{1/\gamma},
  \eea
  where $1/p+1/q+1/q^\prime=1$ and $1/\beta+1/\gamma=1$.
  Fix $\alpha=\beta_0/2$, where $\beta_0$ is defined by
  (\ref{beta}). It follows from $(i)$ of theorem 2.3 that there
  exists some constant $\tau$ depending only on $\alpha$, $n$, $K$
  and $i_0$ such that
  \be\label{i0}\Lambda_\alpha:=\sup_{\|u\|_{1,\tau}\leq 1}\int_M
  \zeta\le(n,\alpha|u|^{\f{n}{n-1}}\ri)dv_g<+\infty.\ee
  Let $r$ be a positive constant to be determined later. Now suppose $\|u\|_E=r$.
  It is easy to see that $\|u\|_{1,\tau}\leq r+\tau r/v_0^{1/n}$. Clearly one can
  select $r$ sufficiently small such that
  $q^\prime\alpha_0\|u\|_{1,\tau}^{n/(n-1)}<\alpha$ and
  $\gamma\alpha_0\|u\|_{1,\tau}^{n/(n-1)}<\alpha$. It follows from (\ref{i0}) that
   $$
   \sup_{\|u\|_E=r}\int_M\zeta\le(n,q^\prime\alpha_0|u|^{\f{n}{n-1}}\ri)dv_g\leq
   \Lambda_\alpha$$
   and
   $$\sup_{\|u\|_E=r}\int_M\zeta\le(n,\gamma\alpha_0|u|^{\f{n}{n-1}}\ri)dv_g\leq
   \Lambda_\alpha,$$
  provided that $r$ is chosen sufficiently small. Inserting these
  two inequalities into (\ref{mmm}), then using the embedding $E\hookrightarrow
  L^s(M)$ for all $s\geq n$ (proposition 7.1) and (\ref{lam-p}), we obtain
  $$J(u)\geq \f{\theta}{n}\|u\|_E^n-\tilde{C}\|u\|_E^{n+1}$$
  for some constant $\tilde{C}$ depending only on $\alpha$, $n$, $K$ and $i_0$,
  provided that $\|u\|_E$ is sufficiently small.
  This gives the desired result. $\hfill\Box$\\

  To estimate the min-max level of $J$, we state the following:\\

  \noindent{\bf Lemma 7.5.} {\it Assume $(f_5)$. There exists some nonnegative function
  $u^*\in E$ such that
  $$\sup_{t\geq 0}J(tu^*)<\f{1}{n}\le(\f{(p-1)\alpha_n}{p\alpha_0}\ri)^{n-1}.$$}

  \noindent{\it Proof.} Let $u^*$ be given by proposition 7.2, namely $u^*\geq
  0$,  $\|u^*\|_E=S_q$, and $\int_M\phi |u^*|^qdv_g=1$. Then
 for any $t\geq 0$ there holds
 \bna
  J(tu^*)&=&\f{1}{n}\|tu^*\|_E^n-\int_{M}\phi(x) F(x,tu^*)dv_g\\
  &\leq& \f{S_q^n}{n}t^n-\f{C_q}{q}t^q\\
  &\leq&\f{q-n}{n q}\f{S_q^{{nq}/{(q-n)}}}{C_q^{{n}/{(q-n)}}}\\
  &<&\f{1}{n}\le(\f{(p-1)\alpha_n}{p\alpha_0}\ri)^{n-1}.
 \ena
 Here we have used the hypothesis $(f_5)$. $\hfill\Box$\\

  Adapting the proof of (\cite{Yang}, lemma 3.4), we obtain the
  following compactness result.\\

  \noindent{\bf Lemma 7.6.} {\it Assume  $(f_1)$, $(f_2)$ and $(f_3)$.
  Let $(u_j)\subset E$ be an arbitrary Palais-Smale sequence of $J$, i.e.,
   \be\label{PS}J(u_j)\ra c,\,\, J^\prime(u_j)\ra 0 \,\, {\rm in}\,\,
   E^*
   \,\,{\rm as}\,\, j\ra \infty,\ee
   where $E^*$ denotes the dual space of $E$.
   Then there exist a subsequence of $(u_j)$ (still denoted by $(u_j)$) and $u\in E$ such that
   $u_j\rightharpoonup u$ weakly in $E$, $u_j\ra u$ strongly in
   $L^q(M)$ for all $q\geq n$, and
   \bna\le\{\begin{array}{lll}
   \nabla u_j(x)\ra\nabla u(x)\quad{\rm
   a.\,\,e.\,\,\,in}\quad M\\[1.5ex]
    {\phi(x)F(x,\,u_j)}\ra {\phi(x)F(x,\,u)}\,\, {\rm strongly\,\,in}\,\,
   L^1(M).
   \end{array}\ri.\ena
   Furthermore $u$ is a weak solution of (\ref{equa}). }\\

   \noindent{\it Proof.}
    Assume $(u_j)$ is a Palais-Smale sequence of $J$. By $(\ref{PS})$, we have
   \bea\label{1}
   &&\f{1}{n}\|u_j\|_E^n-
   \int_{M}\phi(x)F(x,u_j)dv_g\ra
   c\,\,{\rm as}\,\,j\ra\infty,\\\label{2}
   &&\le|\int_{M}\le(|\nabla_g u_j|^{n-2}\nabla_g
   u_j\nabla_g\psi+v|u_j|^{n-2}u_j\psi\ri)dv_g
   -\int_{M}\phi(x)f(x,u_j)\psi dv_g\ri|\leq \sigma_j\|\psi\|_E\qquad\quad
   \eea
   for all $\psi\in E$, where $\sigma_j\ra 0$ as $j\ra\infty$. Note that $f(x,s)\equiv 0$ for all
   $(x,s)\in M\times(-\infty,0]$. By $(f_2)$, we have $0\leq \mu F(x,u_j)\leq
   u_jf(x,u_j)$ for some $\mu>n$. Taking $\psi=u_j$ in (\ref{2}) and multiplying
   (\ref{1}) by $\mu$, we have
   \bna
   \le(\f{\mu}{n}-1\ri)\|u_j\|_E^n&\leq&\mu |c|+\int_{M}\phi(x) \le(\mu F(x,u_j)-f(x,u_j)
   u_j\ri)dv_g+\sigma_j\|u_j\|_E\\
   &\leq&\mu |c|+\sigma_j\|u_j\|_E.
   \ena
   Therefore $\|u_j\|_E$ is bounded. It then follows from  (\ref{1}) and (\ref{2}) that
   \be\label{c}\int_{M}\phi(x)f(x,u_j)u_jdv_g\leq C,\quad \int_{M}\phi(x)F(x,u_j)dv_g\leq
   C\ee for some constant $C$ depending only on $\mu$, $n$ and $c$.
   By proposition 7.1, there exists some $u\in E$ such that $u_j\rightharpoonup
   u$ weakly in $E$, $u_j\ra u$ strongly in $L^q(M)$ for any $q\geq
   n$, and $u_j\ra u$ almost everywhere in $M$. By $(f_3)$,
   there exist positive constants $A_1$ and $R_1$ such that $F(x,s)\leq A_1 f(x,s)$ for all
   $s\geq R_1$. Particularly for any $A>R_1$ there holds
   \be\label{F-g}F(x,s)\leq A_1f(x,s),\quad\forall s\geq A.\ee
   Now we prove that $\phi(x)F(x,u_j)\ra \phi F(x,u)$ strongly in $L^1(M)$. To this end, for any $\epsilon>0$,
   we take $A>\max\{A_1C/\epsilon,R_1\}$, where $C$ is given by (\ref{c}). Then we have by (\ref{F-g})
   \be\label{g-A}
   \int_{|u_j|>A}\phi(x)F(x,u_j)dv_g\leq
   \f{A_1}{A}\int_M\phi(x)f(x,u_j)u_jdv_g<\epsilon.
   \ee
   In the same way
   \be\label{g-A-u}\int_{|u|>A}\phi(x)F(x,u)dv_g<\epsilon.\ee
   By $(f_1)$, we have for $(x,s)\in M\times [0,\infty)$
   \bna
    f(x,s)&\leq&
    b_1s^{n-1}+b_2\zeta\le(n,\alpha_0s^{\f{n}{n-1}}\ri)\\
    &=&b_1s^{n-1}+b_2s^{n}\sum_{k=n-1}^\infty\f{\alpha_0^k
    s^{\f{n}{n-1}(k-n+1)}}{k!}\\
    &\leq&b_1s^{n-1}+b_2s^{n}\alpha_0^{n-1}e^{\alpha_0s^{\f{n}{n-1}}}.
   \ena
   Hence for all $(x,s)\in M\times[0,A]$ there holds
   $$f(x,s)\leq \le(b_1+b_2\alpha_0^{n-1}Ae^{\alpha_0A^{\f{n}{n-1}}}\ri)s^{n-1}.$$
   It follows that
   $$F(x,s)\leq \f{b_1+b_2\alpha_0^{n-1}Ae^{\alpha_0A^{\f{n}{n-1}}}}{n}s^{n},\quad\forall
   s\in[0,A].$$
   for all $(x,s)\in M\times[0,A]$, which implies
   \be\label{le}|\phi(x)\chi_{\{|u_j|\leq A\}}(x)F(x,u_j)|\leq C_1\phi(x)|u_j|^n,\ee
   where $C_1={(b_1+b_2\alpha_0^{n-1}Ae^{\alpha_0A^{{n}/{(n-1)}}})}/{n}$ and
   $\chi_{\{|u_j|\leq A\}}(x)$ denotes the characteristic
   function of the set $\{x\in M: |u_j(x)|\leq A\}$.  By an
   inequality $||a|^n-|b|^n|\leq n|a-b|(|a|^{n-1}+|b|^{n-1})\,(\forall
   a,b\in\mathbb{R})$ and H\"older's inequality, we get
   \bna
    \int_M\phi||u_j|^n-|u|^n|dv_g&\leq&n\int_M\phi|u_j-u|(|u_j|^{n-1}+|u|^{n-1})dv_g\\
    &\leq&n\le(\int_M\phi|u_j-u|^ndv_g\ri)^{\f{1}{n}}\le\{\le(\int_M\phi|u_j|^ndv_g\ri)^{1-\f{1}{n}}
    +\le(\int_M\phi|u|^ndv_g\ri)^{1-\f{1}{n}}\ri\}.
   \ena
   Hence $\phi|u_j|^n\ra \phi|u|^n$ in $L^1(M)$ since $u_j\ra u$
   strongly in $L^n(M)$. In view of (\ref{le}), we conclude from the generalized
   Lebesgue's dominated convergence theorem
   $$\lim_{j\ra\infty}\int_M\phi(x)\chi_{\{|u_j|\leq A\}}(x)F(x,u_j)dv_g
   =\int_M\phi(x)\chi_{\{|u|\leq A\}}(x)F(x,u)dv_g.$$
   This together with (\ref{g-A}) and (\ref{g-A-u}) implies that
   there exists some $m\in \mathbb{N}$ such that when $j>m$ there
   holds
   $$\le|\int_M\phi F(x,u_j)dv_g-\int_M\phi F(x,u)dv_g\ri|<3\epsilon.$$
   Therefore
   $$\lim_{j\ra\infty}\int_M\phi F(x,u_j)dv_g=\int_M\phi F(x,u)dv_g.$$
   Using the same method as that of proving (\cite{Adi-Yang}, (4.26)), we have
   $\nabla_g u_j(x)\ra\nabla_g u(x)$ for almost every $x\in M$
   and
   $$|\nabla_g u_j|^{n-2}\nabla_g u_j\rightharpoonup |\nabla_g u|^{n-2}\nabla_g u
   \quad{\rm weakly\,\,in}\quad \le(L^{\f{n}{n-1}}
   (M)\ri)^n.$$
   Passing to the limit $j\ra \infty$ in $(\ref{2})$, we obtain
   $$\int_{M}\le(|\nabla_g u|^{n-2}\nabla_g u\nabla\psi+v|u|^{n-2}u\psi\ri)dv_g-
   \int_{M}\phi(x)f(x,u)\psi dv_g=0$$
   for all $\psi\in C_0^\infty(M)$. Since $C_0^\infty(M)$ is dense in $E$ under the norm
   $\|\cdot\|_E$, $u$ is a weak solution of
   $(\ref{equa})$.  $\hfill\Box$\\

   We say more words on lemma 7.6. Suppose $(M,g)$ is the standard euclidian space $\mathbb{R}^n$ and
   $\phi(x)=|x|^{-\beta}$, $0\leq \beta<n$. The author \cite{Yang} proved that
   $\phi F(x,u_j)\ra \phi F(x,u)$ in $L^1(\mathbb{R}^n)$ under the assumption
   $E\hookrightarrow L^q(\mathbb{R}^n)$ compactly for all $q\geq 1$. While Lam-Lu \cite{LamLu} observed that
   the convergence still holds  under the assumption $E\hookrightarrow
   L^q(\mathbb{R}^n)$ for all $q\geq n$. Here we generalized these two situations.  \\

   The following lemma is a nontrivial consequence of
   theorem 2.3. It is sufficient for our use when we consider the existence and multiplicity results
   for problems (\ref{equa}) and (\ref{equa1}).\\

  \noindent{\bf Lemma 7.7.} {\it Let $(u_j)\subset E$ be any sequence of functions
  satisfying
  $\|u_j\|_E\leq 1$, $u_j\rightharpoonup u_0$ weakly in $E$, $\nabla_gu_j\ra
  \nabla_gu_0$ almost everywhere in M, and $u_j\ra u_0$ strongly in $L^n(M)$
  as $j\ra\infty$. Then\\
  $(i)$ for any $\alpha: 0<\alpha<\alpha_n$, there holds
 \be\label{1111}\sup_j\int_{M}\zeta\le(n,\alpha|u_j|^{\f{n}{n-1}}\ri)dv_g<\infty;\ee
 $(ii)$ for any
 $\alpha:0<\alpha<\alpha_n$ and $q: 0<q<(1-\|u_0\|_E^n)^{-1/(n-1)}$, there holds
 \be\label{co}\sup_j\int_M\zeta\le(n,q\alpha|u_j|^{\f{n}{n-1}}\ri)dv_g<\infty.\ee}

 \noindent {\it Proof.} $(i)$ For any fixed $\alpha: 0<\alpha<\alpha_n$, it follows
 from part $(i)$ of theorem 2.3 that there exists a positive constant $\tau_{\alpha}$
 depending only on $\alpha$,
 $n$, $K$ and $i_0$ such that
 \be\label{thm}\mathcal{B}_\alpha=\sup_{u\in W^{1,n}(M),\,\|u\|_{1,\tau_{\alpha}}\leq
 1}\int_M\zeta\le(n,\alpha|u|^{\f{n}{n-1}}\ri)dv_g<\infty.\ee
 Note that $v\geq v_0$ in $M$. Since $\|u_j\|_E\leq 1$, we get
 \bna
  \|u_j\|_{1,\tau_{\alpha}}=\le(\int_M|\nabla_gu_j|^ndv_g\ri)^{\f{1}{n}}+\tau_{\alpha}
  \le(\int_M|u_j|^ndv_g\ri)^{\f{1}{n}}
  \leq 1+\f{\tau_{\alpha}}{v_0^{1/n}}.
  \ena
  There exists some small positive number $\alpha_0$ such that $\alpha_0\|u_j\|_{1,\tau_{\alpha}}^{\f{n}{n-1}}\leq
  \alpha$. Hence by (\ref{thm}), there holds
   $$\sup_j\int_M\zeta\le(n,\alpha_0|u_j|^{\f{n}{n-1}}\ri)dv_g\leq
   \sup_{j}\int_M\zeta\le(n,\alpha\le|\f{u_j}{\|u_j\|_{1,\tau_\alpha}}\ri|^{\f{n}{n-1}}\ri)dv_g\leq
   \mathcal{B}_\alpha.$$
  This allows us to define
  $$\alpha^*=\sup\le\{\alpha: \sup_j\int_M\zeta\le(n,\alpha|u_j|^{\f{n}{n-1}}\ri)dv_g<\infty\ri\}.$$
  To prove (\ref{1111}), it suffices to prove that $\alpha^*\geq \alpha_n$. Suppose not, we have
  $\alpha^*<\alpha_n$.  Take two constants $\alpha^\prime$ and $\alpha^{\prime\prime}$ such that
  $\alpha^*<\alpha^\prime<\alpha^{\prime\prime}<\alpha_n$. By part $(i)$ of theorem 2.3 again, there exists some constant
  $\tau_{\alpha^{\prime\prime}}$ depending only on
  $\alpha^{\prime\prime}$, $n$, $K$ and $i_0$ such that
  \be\label{thm-1}\mathcal{B}_{\alpha^{\prime\prime}}=\sup_{u\in W^{1,n}(M),\,\|u\|_{1,\tau_{\alpha^{\prime\prime}}}\leq
 1}\int_M\zeta\le(n,\alpha^{\prime\prime}|u|^{\f{n}{n-1}}\ri)dv_g<\infty.\ee
  Since $u_j\ra u_0$ strongly in $L^n(M)$ and $\nabla_gu_j\ra\nabla_gu_0$ a. e. in M, we
  obtain by using Brezis-Lieb's lemma \cite{BL}
  $$\|u_j-u_0\|_{1,\tau_{\alpha\,^{\prime\prime}}}=\le(\int_M|\nabla_gu_j|^ndv_g-\int_M|\nabla_gu_0|^ndv_g\ri)^{{1}/{n}}
  +o_j(1),$$ where $o_j(1)\ra 0$ as $j\ra\infty$. Since $u_j\rightharpoonup
  u_0$ weakly in $E$, there holds
  $$\lim_{j\ra+\infty}\int_M|\nabla_g u_0|^{n-2}\nabla_gu_0\nabla_gu_j\,dv_g=\int_M|\nabla_gu_0|^ndv_g.$$
  This immediately implies that
  $$\int_M|\nabla_gu_0|^ndv_g\leq \limsup_{j\ra+\infty}\int_M|\nabla_gu_j|^ndv_g\leq 1.$$
  Hence
  $$\|u_j-u_0\|_{1,\tau_{\alpha\,^{\prime\prime}}}\leq 1+o_j(1).$$
  It follows from (\ref{Y2}) that for any $\epsilon>0$ there
  exists some constant $\tilde{c}$ depending only on $\epsilon$ and $n$
  such that
  \bea\label{nu}
   \zeta\le(n,\alpha^\prime|u_j|^{\f{n}{n-1}}\ri)\leq\f{1}{\mu}\zeta\le(n,\alpha^\prime(1+\epsilon)\mu|u_j-u_0|^{\f{n}{n-1}}\ri)+
   \f{1}{\nu}\zeta\le(n,\alpha^\prime\tilde{c}\nu|u_0|^{\f{n}{n-1}}\ri),
  \eea
  where $1/\mu+1/\nu=1$.
  Choosing $\epsilon$ sufficiently small and $\mu$ sufficiently
  close to $1$ such that
  $$\alpha^\prime(1+\epsilon)\mu\|u_j-u_0\|_{1,\tau_{\alpha^{\prime\prime}}}^{\f{n}{n-1}}<
  \alpha^{\prime\prime},$$
  provided that $j$ is sufficiently large. This together with (\ref{thm-1}) implies that
  \be\label{tm3}\sup_j\int_M\zeta\le(n,\alpha^\prime(1+\epsilon)\mu|u_j-u_0|^{\f{n}{n-1}}\ri)dv_g
  \leq \mathcal{B}_{\alpha^{\prime\prime}}.\ee
  In addition, we have by part $(iii)$ of theorem 2.3 that
  \be\label{tm4}\int_M\zeta\le(n,\alpha^\prime\tilde{c}\nu|u_0|^{\f{n}{n-1}}\ri)dv_g<+\infty.\ee
  Inserting (\ref{tm3}) and (\ref{tm4}) into (\ref{nu}), we get
  $$\sup_j\int_M\zeta\le(n,\alpha^\prime|u_j|^{\f{n}{n-1}}\ri)dv_g<+\infty,$$
  which contradicts the definition of $\alpha^*$ and thus ends the proof
  of part $(i)$.  \\

  \noindent $(ii)$ Given any $\alpha:0<\alpha<\alpha_n$ and any $q:0<q<(1-\|u_0\|_E^n)^{-1/(n-1)}$.
   By (\ref{Y2}), $\forall\epsilon>0$, there exist constants $\tilde{c}>0$, $\mu>1$
   and $\nu>1$ $(1/\mu+1/\nu=1)$ such that
  \bna
   \int_M\zeta\le(n,q\alpha|u_j|^{\f{n}{n-1}}\ri)dv_g\leq\f{1}{\mu}
   \int_M\zeta\le(n,q\alpha(1+\epsilon)\mu|u_j-u_0|^{\f{n}{n-1}}\ri)dv_g+
   \f{1}{\nu}\int_M\zeta\le(n,q\alpha\tilde{c}\nu|u_0|^{\f{n}{n-1}}\ri)dv_g.
  \ena
  By Brezis-Lieb's lemma \cite{BL},
  $$\|u_j-u_0\|_E^{\f{n}{n-1}}\leq (1-\|u_0\|_E^n)^{\f{1}{n-1}}+o_j(1).$$
  If we choose $\epsilon$
  sufficiently small and $\mu$ sufficiently close to $1$ such that
  $$q\alpha(1+\epsilon)\mu\|u_j-u_0\|_E^{\f{n}{n-1}}\leq(\alpha+\alpha_n)/2,$$
  provided that $j$ is sufficiently large. It then follows from part $(i)$
  that
  $$\sup_j\int_M\zeta\le(n,q\alpha(1+\epsilon)\mu|u_j-u_0|^{\f{n}{n-1}}\ri)dv_g<+\infty.$$
  By part $(iii)$ of theorem 2.3, we have
  $$\int_M\zeta\le(n,q\alpha\tilde{c}\nu|u_0|^{\f{n}{n-1}}\ri)dv_g<+\infty.$$
  Therefore (\ref{co}) holds. $\hfill\Box$\\

 {\bf Remark 7.8.} In lemma 7.7, if $u_0\equiv0$, then the conclusion of $(ii)$ is weaker than
 that of $(i)$. If $u_0\not\equiv 0$, then the conclusion of $(i)$ is a special case of
 that of $(ii)$. If $(M,g)$ has dimension two, the assumption
 $\nabla_g u_j\ra \nabla_g u_0$ almost everywhere in $M$ can be removed.\\

 {\it Proof of theorem 2.7.}
 It follows from lemma 7.3 and lemma 7.4 that $J$ satisfies all the hypothesis of
 the mountain-pass theorem except for the Palais-Smale condition: $J\in
   \mathcal{C}^1(E,\mathbb{R})$; $J(0)=0$; $J(u)\geq \delta>0$ when
  $\|u\|_E=r$; $J(e)<0$ for some  $e\in E$  with
   $\|e\|_E>r$.
   Then using the mountain-pass theorem
   without the Palais-Smale condition \cite{Rabin}, we can find a sequence
   $(u_j)$ in $E$ such that
   $$J(u_j)\ra c>0,\quad J^\prime(u_j)\ra 0\,\,{\rm in}\,\, E^*,$$
   where
   $$c=\min_{\gamma\in\Gamma}\max_{u\in\gamma}J(u)\geq \delta$$
   is the min-max value of $J$, where $\Gamma=\{\gamma\in\mathcal{C}([0,1],E): \gamma(0)=0,
   \gamma(1)=e\}$. This is equivalent to (\ref{1}) and $(\ref{2})$. By lemma 7.6,
   up to a subsequence, there holds
   \be\label{Fcon}\le\{\begin{array}{lll}
   u_j\rightharpoonup u \,\,{\rm weakly\,\, in}\,\, E\\[1.5ex]  u_j\ra
   u \,\,{\rm strongly\,\, in}\,\, L^q(M),\,\,\forall q\geq
   n\\[1.5ex]
   \lim\limits_{j\ra\infty}\int_{M}
   \phi(x){F(x,u_j)}dv_g=\int_{M}\phi(x){F(x,u)}dv_g\\
   [1.5ex] u\,\,{\rm is\,\, a\,\,weak\,\,solution\,\,of}\,\,(\ref{equa}).\end{array}
   \ri.\ee
   Now suppose by contradiction $u\equiv0$. Since $F(x,0)= 0$ for all $x\in M$,  it follows from
   (\ref{1}) and (\ref{Fcon}) that
   \be\label{N}\lim_{j\ra\infty}\|u_j\|_E^n=nc>0.\ee
   By lemma 7.5, $0<c<\f{1}{n}\le(\f{(p-1)\alpha_n}{p\alpha_0}\ri)^{n-1}$.
   Thus there exists some $\eta_0>0$ and $m>0$ such that
   $\|u_j\|_E^n\leq \le(\f{p-1}{p}\f{\alpha_n}{\alpha_0}-\eta_0\ri)^{n-1}$ for all $j>m$.
    Choose $q>1$ sufficiently close to $1$ such
   that $q\alpha_0\|u_j\|_E^{\f{n}{n-1}}\leq (1-1/p)\alpha_n-\alpha_0\eta_0/2$ for all
   $j>m$.    By $(f_1)$,
   $$|f(x,u_j)u_j|\leq b_1|u_j|^n+b_2|u_j|\zeta\le(n,\alpha_0|u_j|^{\f{n}{n-1}}\ri).$$
    It follows from (\ref{Y1}), H\"older's inequality, and
    part $(i)$ of lemma 7.7 that
   \bna
   \int_{M}\phi{|f(x,u_j)u_j|}dv_g&\leq&b_1\int_{M}
   \phi|u_j|^ndv_g+b_2\int_{M}\phi|u_j|\zeta\le(n,\alpha_0|u_j|^{\f{n}{n-1}}\ri)dv_g\\
   &\leq&b_1\int_{M}\phi{|u_j|^n}dv_g+b_2\le(\int_{M}\phi|u_j|^{q^\prime}
   dv_g\ri)^{1/{q^\prime}}
   \le(\int_{M}\phi\zeta\le(n,q\alpha_0|u_j|^{\f{n}{n-1}}\ri)dv_g\ri)^{1/{q}}\\
   &\leq&b_1\int_{M}\phi{|u_j|^n}dv_g+C\le(\int_{M}\phi{|u_j|^{q^\prime}}dv_g\ri)^{1/{q^\prime}}
   \ra 0\quad{\rm as}\quad j\ra\infty,
   \ena
   where $1/q+1/q^\prime=1$ and $C$ is some constant which is
   independent of $j$.
   Here we have used (\ref{Fcon}) again (precisely $u_j\ra u$ in
   $L^r(\mathbb{R}^N)$ for all $r\geq n$) in the above estimates.
    Inserting this into (\ref{2}) with $\psi=u_j$, we have
   $$\|u_j\|_E\ra 0\quad{\rm as}\quad j\ra\infty,$$
   which contradicts (\ref{N}). Therefore $u\not\equiv 0$ and we obtain a nontrivial
   weak solution of (\ref{equa}).
   Finally $u$ is nonnegative since
   $f(x,s)\equiv 0$ for all $(x,s)\in M\times(-\infty,0]$. $\hfill\Box$\\

   {\it Proof of theorem 2.10.} Since the proof is very similar to that of (\cite{Yang}, theorem 1.2),
   we only give its sketch and emphasize the difference between these two situations.
   Instead of $J:E\ra\mathbb{R}$ defined by (\ref{J}), we consider
   functionals for all $u\in E$ and $\epsilon>0$
   $$J_{\epsilon}(u)=\f{1}{n}\|u\|_E^n-\int_{M}\phi(x){F(x,u)}dv_g-\epsilon\int_{M}hudv_g.$$
   Firstly, lemma 7.6 still holds if we replace $J$ by $J_\epsilon$. Namely
   for any Palais-Smale sequence $(u_j)\subset E$ of $J_\epsilon$,
    there exist a subsequence of $(u_j)$ (still denoted by $(u_j)$) and $u\in E$ such that
   $u_j\rightharpoonup u$ weakly in $E$, $u_j\ra u$ strongly in
   $L^q(M)$ for all $q\geq n$, and
   \be\label{ps1}\le\{\begin{array}{lll}
   \nabla_g u_j(x)\ra\nabla_g u(x)\quad{\rm
   a.\,\,e.\,\,\,in}\quad M\\[1.5ex]
    {\phi(x)F(x,\,u_j)}\ra {\phi(x)F(x,\,u)}\,\, {\rm strongly\,\,in}\,\,
   L^1(M)\\[1.5ex]
   u\,\,{\rm is\,\,a\,\,weak\,\,solution\,\,of }\,\,(\ref{equa1}).
   \end{array}\ri.\ee
   Secondly, using the same method in the first two steps of the proof of (\cite{Yang},
   theorem 1.2), we have the following:\\

   \noindent $(a)$ there exist constants $\epsilon_1>0$,
   $\delta>0$, and a sequence of functions $(v_j)\subset E$ such
   that $J_\epsilon(v_j)\ra c_M$ and $J^\prime_\epsilon(v_j)\ra 0$ as
   $j\ra\infty$, provided that $0<\epsilon<\epsilon_1$. In addition, $v_j$ is bounded in $E$,
   $v_j\rightharpoonup u_M$ weakly in $E$
   and $u_M$ is a weak solution of (\ref{equa1}). Here $c_M$ is the
   min-max value of $J_\epsilon$ and
   satisfies
   \be\label{cM1}0<c_M<\f{1}{n}\le(1-\f{1}{p}\ri)^{n-1}\le(\f{\alpha_n}{\alpha_0}\ri)^{n-1}-\delta;\ee
   $(b)$ there exists a constant $\epsilon_2: 0<\epsilon_2<\epsilon_1$
   such that for any $\epsilon: 0<\epsilon<\epsilon_2$, there exist positive constant
   $r_\epsilon$ with $r_\epsilon\ra 0$ as $\epsilon\ra 0$ and
   sequence $(u_j)\subset E$ such that
   $$J_\epsilon(u_j)\ra c_\epsilon=\inf_{\|u\|_E\leq r_\epsilon}J_{\epsilon}(u)<0,\quad
   J_{\epsilon}^\prime(u_j)\ra 0\quad{\rm in}\quad E^*\quad{\rm as}\quad j\ra\infty.$$
   In addition, $u_j\ra u_0$ strongly in $E$, where
   $u_0$ is a
weak solution of (\ref{equa1}) with
$J_{\epsilon}(u_0)=c_\epsilon$.\\

\noindent Thirdly, there exists $\epsilon_0:
0<\epsilon_0<\epsilon_2$ such that if $0<\epsilon<\epsilon_0$,
 then $u_M\not\equiv u_0$. Suppose by contradiction that $u_M\equiv u_0$.
 Then $v_j\rightharpoonup u_0$ weakly in $E$. By $(a)$,
 \be\label{mt} J_{\epsilon}(v_j)\ra c_M>0,\quad
  |\langle  J_{\epsilon}^\prime(v_j),\varphi\rangle|\leq \gamma_j\|\varphi\|_E\ee
  with $\gamma_j\ra 0$ as $j\ra\infty$. On one hand we have by (\ref{ps1}),
 \be\label{cv}\int_{M}\phi(x)F(x,v_j)dv_g\ra\int_{M}\phi(x)F(x,u_0)dv_g
 \quad{\rm as}\quad j\ra\infty.\ee
  On the other hand,
 since $v_j\rightharpoonup u_0$ weakly in $E$ and $h\in E^*$, it follows that
 \be\label{ch}
 \int_{M}hv_jdv_g\ra \int_{M}hu_0dv_g\quad{\rm as}\quad j\ra\infty.
 \ee
 Inserting (\ref{cv}) and (\ref{ch}) into the first equality of (\ref{mt}), we obtain
 \be\label{cM}
 \f{1}{n}\|v_j\|_E^n=c_M+\int_{M}\phi(x)F(x,u_0)dv_g+\epsilon\int_{M}hu_0dv_g+o_j(1).
 \ee
 In the same way, one can derive
 \be\label{c0}
 \f{1}{n}\|u_j\|_E^n=c_\epsilon+\int_{M}\phi(x)F(x,u_0)dv_g+
 \epsilon\int_{M}hu_0dv_g+o_j(1).
 \ee
 Combining (\ref{cM}) and (\ref{c0}), we have
 \be\label{Lp}
 \|v_j\|_E^n-\|u_0\|_E^n=n\le(c_M-c_\epsilon+o_j(1)\ri).
 \ee
 From $(b)$, we know that $c_\epsilon\ra 0$ as $\epsilon\ra 0$. This together with (\ref{cM1}) leads to the existence
 of $\epsilon_0: 0<\epsilon_0<\epsilon_2$ such that if $0<\epsilon<\epsilon_0$, then
 \be\label{diff}0<c_M-c_\epsilon<\f{1}{n}\le(\f{p-1}{p}\f{\alpha_n}{\alpha_0}\ri)^{n-1}.\ee
  Write
 $$w_j=\f{v_j}{\|v_j\|_E},\quad w_0=\f{u_0}{\le(\|u_0\|_E^n+n(c_M-c_\epsilon)\ri)^{1/n}}.$$
 It follows from (\ref{Lp}) and $v_j\rightharpoonup u_0$ weakly in
 $E$ that
 $w_j\rightharpoonup w_0$ weakly in $E$. Note that
 $$\int_{M}\phi(x)\zeta\le(n,\alpha_0|v_j|^{n/(n-1)}\ri)dv_g=\int_{M}
 \phi(x)\zeta\le(n,\alpha_0\|v_j\|_E^{{n}/{(n-1)}}|w_j|^{n/(n-1)}\ri)dv_g.$$
 By (\ref{Lp}) and (\ref{diff}), a straightforward calculation shows
 $$\lim_{j\ra\infty}\alpha_0\|v_j\|_E^{\f{n}{n-1}}\le(1-\|w_0\|_E^n\ri)^{\f{1}{n-1}}<
 \le(1-\f{1}{p}\ri)\alpha_n.$$
 Hence lemma 7.7 together with (\ref{Y2}) implies that
 $\phi(x)\zeta\le(n,\alpha_0|v_j|^{n/(n-1)}\ri)$
  is bounded in $L^q(M)$ for some
 $q: 1<q<n/(n-1)$. By $(f_1)$,
 $$|f(x,v_j)|\leq b_1|v_j|^{n-1}+b_2\zeta(n,\alpha_0|v_j|^{\f{n}{n-1}}).$$
 By the definition of $\zeta$ there exists a constant $c>0$ such
 that
 $$
 |f(x,v_j)\chi_{\{|v_j|\leq 1\}}(x)|\leq c|v_j|^{n-1},\quad
 |f(x,v_j)\chi_{\{|v_j|> 1\}}(x)|\leq
 c\zeta(n,\alpha_0|v_j|^{\f{n}{n-1}}),
 $$
 where $\chi_{B}$ denotes the characteristic function of $B\subset M$.
 Hence
 \bna\le|\int_M\phi(x) f(x,v_j)(v_j-u_0)dv_g\ri|&\leq& c\int_M\phi(x)\le(|v_j|^{n-1}+
 \zeta(n,\alpha_0|v_j|^{\f{n}{n-1}})\ri)|v_j-u_0|dv_g\\
 &\leq& c\le\|\phi|v_j|^{n-1}\ri\|_{L^{\f{n}{n-1}}(M)}
  \|v_j-u_0\|_{L^n(M)}\\
  &&+c\le\|\phi
 \zeta(n,\alpha_0|v_j|^{\f{n}{n-1}})\ri\|_{L^q(M)}\|v_j-u_0\|_{L^{q^\prime}(M)}.\ena
 Since $1<q<n/(n-1)$, we have $q^\prime>n$.
  Then it follows from the compact embedding $E\hookrightarrow L^r(M)$
 for all $r\geq n$ that
 \be\label{f-c}
 \lim_{j\ra\infty}\int_{M}\phi(x)f(x,v_j)(v_j-u_0)dv_g= 0.
 \ee
 Taking $\varphi=v_j-u_0$ in (\ref{mt}), we have by using (\ref{ch}) and (\ref{f-c}) that
 \be\label{cc}
  \int_{M}\le(|\nabla_g v_j|^{n-2}\nabla_g
  v_j\nabla_g(v_j-u_0)+v(x)|v_j|^{n-2}v_j(v_j-u_0)\ri)dv_g
  \ra 0.
 \ee
 However the fact $v_n\rightharpoonup u_0$ weakly in $E$ leads to
 \be\label{cc1}
  \int_{M}\le(|\nabla_g u_0|^{n-2}\nabla_g u_0\nabla_g(v_j-u_0)+
  v(x)|u_0|^{n-2}u_0(v_j-u_0)\ri)dv_g
  \ra 0.
 \ee
 Subtracting (\ref{cc1}) from (\ref{cc}), using the well known inequality (see \cite{Lind}, chapter 10)
  $$2^{n-1}|b-a|^n\leq
 \langle|b|^{n-2}b-|a|^{n-2}a, b-a\rangle,\quad\forall a,b\in\mathbb{R}^n,$$ we have
   $\|v_j-u_0\|_E^n\ra
 0$ as $j\ra\infty$. This together with (\ref{Lp}) implies that
 $c_M=c_\epsilon$,
 which is absurd since $c_M>0$ and $c_\epsilon<0$. Therefore $u_M\not\equiv u_0$. Since
 $f(x,s)\equiv 0$ for all $(x,s)\in M\times(-\infty,0]$, one can easily see that $u_M\geq 0$
 and $u_0\geq 0$. This completes the proof of the theorem.
 $\hfill\Box$\\

 Finally we shall construct examples of $f$'s to
 show that $(f_1)$-$(f_5)$ do not imply $(H_5)$. \\

   \noindent{\it Proof of proposition 2.9.} Let $\phi$ satisfies
   the hypotheses $(\phi_1)$ and $(\phi_2)$, $p>1$ be given in $(\phi_1)$, $l$ be an integer satisfying $l\geq n$,
   $q=nl/(n-1)+1$ and $S_q$ be defined
   by (\ref{Sq}). In view of lemma 7.2, $S_q$ is attained by some
   nonnegative function $u\in E$. Let $C_q$ be a positive number such that
   $$C_q>\le(\f{q-n}{q}\ri)^{{(q-n)}/{n}}\le(\f{p\alpha_0}{(p-1)\alpha_n}\ri)^{(q-n)(n-1)/n}S_q^q.$$
Let $\chi:[0,\infty)\ra\mathbb{R}$ be a smooth function such that
$0\leq \chi\leq 1$, $\chi\equiv 0$ on $[0,A]$, $\chi\equiv 1$ on
$[2A,\infty)$, and $|\chi\,^\prime|\leq 2/A$, where $A$ is a
positive constant to be determined later. We set \bna
f(t)=\le\{\begin{array}{lll}2^ll!C_q\sum_{k=l}^\infty\f{(t^{\f{n}{n-1}}-\chi(t)t^{\f{1}{n-1}})^k}{k!},
\quad &t\geq 0\\
[1.5ex]0, &t<0.\end{array}\ri.\ena

 Now we check $(f_1)$-$(f_5)$ for appropriate choice of $A$ as follows.\\

 \noindent$(f_1)$: If $A>1$, then $0\leq t^{n/(n-1)}-\chi(t)t^{1/(n-1)}\leq
 t^{n/(n-1)}$ for all $t\geq 0$. Thus
 \bna
  f(t)&=&2^ll!C_q\le(e^{t^{n/(n-1)}-\chi(t)t^{1/(n-1)}}-\sum_{k=0}^{l-1}
  \f{(t^{\f{n}{n-1}}-\chi(t)t^{\f{1}{n-1}})^k}{k!}\ri)\\
  &\leq&2^ll!C_q\le(e^{t^{n/(n-1)}}-\sum_{k=0}^{l-1}
  \f{t^{\f{nk}{n-1}}}{k!}\ri)\\
  &\leq&2^ll!C_q\zeta(n,t^{n/(n-1)})
 \ena
 for all $t\geq 0$. So $(f_1)$ is satisfied when $A>1$.\\
 $(f_2)$: When $t\in[0,A]$, we have
$\chi(t)=0$ and
 \be\label{lq}
 \int_0^tf(t)dt
 =2^ll!C_q\sum_{k=l}^\infty\int_0^t\f{t^{\f{nk}{n-1}}}{k!}dt\leq
 {2^ll!C_qt}\sum_{k=l}^\infty\f{t^{\f{nk}{n-1}}}{k!}= t
 f(t).
\ee When $t\geq A$, we {\it claim} that if $A$ is chosen
sufficiently large, say $A\geq 4^{n-1}$, then
 \be\label{gea}\int_A^t\f{(t^{\f{n}{n-1}}-\chi(t)t^{\f{1}{n-1}})^k}{k!}dt\leq \f{(t^{\f{n}{n-1}}
 -\chi(t)t^{\f{1}{n-1}})^{k+1}}{(k+1)!}-\f{A^{\f{n(k+1)}{n-1}}}{(k+1)!}.\ee
 In fact, if we set
$$\gamma(t)=\int_A^t\f{(t^{\f{n}{n-1}}-\chi(t)t^{\f{1}{n-1}})^k}{k!}dt-\f{(t^{\f{n}{n-1}}
 -\chi(t)t^{\f{1}{n-1}})^{k+1}}{(k+1)!}+\f{A^{\f{n(k+1)}{n-1}}}{(k+1)!},$$
 then $\gamma(A)=0$ and
 $$\gamma\,^\prime(t)=\f{(t^{\f{n}{n-1}}-\chi(t)t^{\f{1}{n-1}})^k}{k!}-
 \f{(t^{\f{n}{n-1}}-\chi(t)t^{\f{1}{n-1}})^k}{k!}\le(\f{n}{n-1}t^{\f{1}{n-1}}-
 \chi\,^\prime(t)t^{\f{1}{n-1}}-\f{1}{n-1}\chi(t)t^{\f{1}{n-1}-1}\ri).$$
 Let $A\geq 4^{n-1}$. Then for $t\in [A,\infty)$ there holds
 \bna
 \f{n}{n-1}t^{\f{1}{n-1}}-
 \chi\,^\prime(t)t^{\f{1}{n-1}}-\f{1}{n-1}\chi(t)t^{\f{1}{n-1}-1}&\geq&\le(\f{n}{n-1}-\f{2}{A}\ri)
 A^{\f{1}{n-1}}-\f{1}{n-1}A^{\f{1}{n-1}-1}\\
 &\geq&4\le(\f{n}{n-1}-\f{2}{4(n-1)}-\f{1}{4(n-1)^2}\ri)\\
 &>&1.
 \ena
  Hence $\gamma\,^\prime(t)\leq
 0$  and thus our claim (\ref{gea}) holds. Note that
 \be\label{0a}\int_0^A\f{t^{\f{nk}{n-1}}}{k!}dt=\f{A^{\f{n(k+1)}{n-1}}}{(k+1)!}\f{k+1}{\f{nk}{n-1}+1}A^{-\f{1}{n-1}}\leq
  \f{A^{\f{n(k+1)}{n-1}}}{(k+1)!}.\ee
  It follows from (\ref{gea}) and (\ref{0a}) that when
   $t\geq A$,
  $$\int_0^t\f{(t^{\f{n}{n-1}}-\chi(t)t^{\f{1}{n-1}})^k}{k!}dt\leq \f{(t^{\f{n}{n-1}}
 -\chi(t)t^{\f{1}{n-1}})^{k+1}}{(k+1)!},$$
 and whence
 \be\label{Fv}\int_0^tf(t)dt\leq f(t)\leq \f{1}{\mu}t f(t)\ee for some $\mu>n$.
 This together with (\ref{lq}) implies that $(f_2)$ holds.\\
 \noindent $(f_3)$: Let $A\geq 4^{n-1}$. In view of (\ref{Fv}), when $t\geq A $,
 $$F(t)=\int_0^tf(t)dt\leq f(t).$$
 Hence $(f_3)$ is satisfied.\\
 \noindent $(f_4)$:
 Since $l>n$, we get $F(t)/t^n\ra 0$ as $t\ra 0+$.
 Hence $(f_4)$ holds. \\
 \noindent $(f_5)$: Note that $t^{n/(n-1)}-t^{1/(n-1)}\geq
 t^{n/(n-1)}/2$ for all $t\geq 2$. Let $A\geq 2$. Then for all $t\geq
 A$ there holds
 $$f(t)\geq 2^ll!C_q\f{(t^{\f{n}{n-1}}-\chi(t)t^{\f{1}{n-1}})^l}{l!}\geq 2^lC_q(t^{\f{n}{n-1}}/2)^l
 =C_qt^{q-1}.$$
 When $t\in [0,A]$, we get
 $$f(t)\geq 2^ll!C_q\f{t^{\f{nl}{n-1}}}{l!}=2^lC_q t^{q-1}.$$
 Hence $(f_5)$ is satisfied. In short, $f(t)$ satisfies
 $(f_1)$-$(f_5)$ if $A\geq 4^{n-1}$.

 Finally we check that $(H_5)$ does not hold. When $t\geq 2A$, we
 have
 $$f(t)=2^ll!C_q\le(e^{t^{\f{n}{n-1}}-t^{\f{1}{n-1}}}-\sum_{k=0}^{l-1}\f{(t^{\f{n}{n-1}}-t^{\f{1}{n-1}})^k}{k!}\ri).$$
 It follows that
 $$\lim_{t\ra+\infty}t f(t)e^{-t^{\f{n}{n-1}}}=0.$$
 Thus $f(t)$ does not satisfy $(H_5)$. $\hfill\Box$\\

{\bf Acknowledgements}. This work was partly supported by the NSFC
11171347
  and the NCET program 2008-2011.

\end{document}